\pdfoutput=1
\documentclass[12pt,oneside]{amsart}
\usepackage{placeins}
\usepackage[margin=1in]{geometry}
\usepackage{mathabx}
\usepackage{latexsym,epsfig,enumitem,multicol,wasysym,amsmath,amsthm,amssymb,amsfonts,amscd,ulem,soul,multirow,rotating,pdflscape,float,dsfont}
\setlist{itemsep=.06125in}
 \newcommand{\FR}{\mathcal{F}\!\mathcal{R}}
 
\usepackage{graphicx}
\usepackage[colorlinks,citecolor=red,pagebackref,hypertexnames=false]{hyperref}
\hypersetup{%
  colorlinks = true,
  citecolor=red,
  linkcolor  = red
}
\restylefloat{table}
\numberwithin{equation}{section}
\theoremstyle{plain}
\newtheorem{theorem}{Theorem}[section]
\newtheorem{lemma}[theorem]{Lemma}
\newtheorem{corollary}[theorem]{Corollary}
\newtheorem{proposition}[theorem]{Proposition}

\theoremstyle{definition}
\newtheorem{definition}[theorem]{Definition}

\theoremstyle{remark}
\newtheorem{remark}[theorem]{Remark}

\def\R{\mathbb R}

\date{\today}      
\author{A. Iosevich, Z. Li, E. Palsson, and A. Yavicoli}
\address{Department of Mathematics, University of Rochester, Rochester, NY, USA}
\email{iosevich@gmail.com}
\address{Department of Mathematics, University of Rochester, Rochester, NY, USA}
\email{zli154@ur.rochester.edu} 
\address{Department of Mathematics, Virgina Tech, Blacksburg, VA, USA}
\email{palsson@vt.edu}
\address{Department of Mathematics, University of British Columbia, Vancouver, BC, Canada}
\email{yavicoli@math.ubc.ca} 

\thanks{A.~I. was supported in part by the National Science Foundation under NSF DMS - 2154232.}

\thanks{A.~Y. was supported in part by the Natural Sciences and Engineering Research Council of Canada, NSERC (GR030571 and GR030540).}

\title[Fourier Ratio and Structure of Measures]{The Fourier Ratio: Uncertainty, Restriction, and Approximation for Compactly Supported Measures}

\begin{document} 

\begin{abstract}
We introduce and systematically study a continuous analog of the Fourier ratio for compactly supported Borel measures. For a measure $\mu$ on $\mathbb{R}^d$ and $f \in L^2(\mu)$, we define the Fourier ratio as
\[
\FR_{\mu, R}(f) = \frac{X_{1,\mu,R}(f)}{X_{2,\mu,R}(f)},
\]
where $X_{p,\mu,R}(f)$ is the $L^p$ norm of a regularized Fourier transform at scale $R$. This quantity, which interpolates between $L^1$ and $L^2$ Fourier information, serves as a fundamental parameter connecting uncertainty principles, Fourier restriction theory, and approximation by trigonometric polynomials.

Our first main contribution is a fractal uncertainty principle (Theorem \ref{theorem:fractalsandwich}) that provides two-sided bounds for $\FR_{\mu,R}(f)$ in terms of the covering numbers of the spatial and frequency supports of $f\mu$. This leads to exact signal recovery results (Theorem \ref{theorem:measurerecovery}) under natural geometric conditions.

Second, we prove that a small Fourier ratio implies efficient low-degree trigonometric approximation: the mollified measure $(f\mu)*\psi_{R^{-1}}$ can be approximated in $L^1$, $L^2$, or $L^\infty$ by a trigonometric polynomial whose degree is bounded explicitly by $\FR_{\mu,R}(f)$ (Theorems \ref{theorem:L2approximation}, \ref{Linfinityapproximation}, and \ref{theorem:L1approximation}). 

Third, we establish a sharp contrast between deterministic and random settings via restriction theory. For the arc-length measure on the circle, we show $\FR_{\sigma,R}(f) \lesssim R^{-1/4}$, while for the Laba--Wang random Cantor measure of dimension $1$, we prove the subpolynomial lower bound $\FR_{\mu,R}(f) \gtrsim R^{-\epsilon}$ for every $\epsilon > 0$ (Theorem \ref{theorem:restrictionFR}). Consequently, $(f\mu)*\psi_{R^{-1}}$ cannot be approximated in $L^2$ with fixed accuracy by trigonometric polynomials of degree $o(R^2)$, in stark contrast with the curved case.

Finally, for convex geometry, we show that the degree of polynomial approximation for the surface measure on a convex body is governed by the upper Minkowski dimension of its set of outward unit normals (Corollary \ref{cor:minkowskikicksass}). 

These results unify discrete and continuous Fourier analysis, demonstrating that the Fourier ratio is a central object linking geometric measure theory, harmonic analysis, and approximation theory.
\end{abstract}

\maketitle

\tableofcontents

\section{Introduction}
\label{section:introduction}

The $L^1$-$L^2$ ratio of the Fourier transform---or Fourier ratio---of a discrete signal $h : \mathbb{Z}_N \to \mathbb{C}$, defined as
\[
\frac{\frac{1}{N}\sum_{m \in {\mathbb Z}_N}|\widehat{h}(m)|}{\big(\frac{1}{N}\sum_{m \in {\mathbb Z}_N}|\widehat{h}(m)|^2\big)^{1/2}},
\]

has emerged as a sensitive measure of pseudorandomness and structure \cite{A2025}. Classical inequalities of Talagrand \cite{Talagrand98} and Bourgain \cite{Bourgain89} show that a large Fourier ratio is typical for functions concentrated on random subsets, while a small ratio forces the signal to be well-approximated by a sparse trigonometric polynomial. A fundamental question is whether this discrete phenomenon reflects a broader principle in continuous harmonic analysis.

\subsection{The Core Problem and Main Definition}
Extending this theory to Euclidean space faces an immediate obstacle: for a general Borel measure $\mu$, the object $\widehat{f\mu}$ may not be integrable, obstructing a direct analogy. The principal contribution of this paper is to resolve this by introducing a scale-dependent regularization of the Fourier transform, yielding a robust continuous Fourier ratio that captures both the geometry of $\mu$ and the analytic properties of $f$.

Let $\mu$ be a compactly supported Borel measure on $\mathbb{R}^d$ and $f \in L^2(\mu)$. Let $\psi$ be a smooth, compactly supported approximation to the identity with $\int\psi = 1$, and define $\psi_{\delta}(x) = \delta^{-d} \psi(x/\delta)$. We study the regularized $L^p$-norms of the Fourier transform at scale $R$:
\begin{equation}
\label{eq:Xdef}
X_{p,\mu,R}(f) := \bigg( R^{-d} \int_{\mathbb R^d} \big| \widehat{(f\mu) * \psi_{R^{-1}}}(\xi) \big|^p d\xi \bigg)^{1/p}, \quad p \in [1,\infty).
\end{equation}
The \textbf{Fourier ratio} of $f$ (with respect to $\mu$ at scale $R$) is then
\begin{equation}
\label{eq:FRdef}
\FR_{\mu,R}(f) := \frac{X_{1,\mu,R}(f)}{X_{2,\mu,R}(f)}.
\end{equation}

This quantity is well-defined for all $f \in L^2(\mu)$ and inherits the key invariance and scaling properties of its discrete ancestor, while coupling intimately with the $R^{-1}$-scale geometry of $\operatorname{spt}(\mu)$.

\subsection{Summary of Main Contributions}
We present four interconnected groups of results that establish the Fourier ratio as a central object in analysis.

\subsubsection{Fractal Uncertainty Principle and Exact Recovery}
We prove that $\FR_{\mu,R}(f)$ is squeezed between geometric quantities determined by the spatial support of $f$ and the frequency concentration of $\widehat{f\mu}$.
\begin{itemize}
    \item \textbf{Geometric Lower Bound (Proposition \ref{prop:lowerbound}):} For any $f$,
    \[
    \FR_{\mu,R}(f) \gtrsim \frac{1}{\sqrt{R^d |E_f^{R^{-1}}|}} \approx \frac{1}{\sqrt{\#\{R^{-1}\text{-balls covering } \operatorname{supp}(f)\}}}.
    \]
    \item \textbf{Upper Bound via $L^1$ Concentration (Theorem \ref{theorem:fractalsandwich}):} If $\widehat{(f\mu)_{R^{-1}}}$ is $L^1$-concentrated on a set $X$, then
    \[
    \FR_{\mu,R}(f) \lesssim \sqrt{\frac{|X \cap B_R|}{R^d}}.
    \]
    \item \textbf{Fractal Uncertainty Principle:} Combining these yields
    \[
    (1-\eta)^2 \lesssim |E_f^{R^{-1}}| \cdot |X \cap B_R|,
    \]
    a continuous analogue of the Bourgain--Dyatlov framework. This leads to \textbf{exact signal recovery} (Theorem \ref{theorem:measurerecovery}) when certain high-frequency data are missing, provided the spatial and frequency supports are sufficiently sparse.
\end{itemize}

\subsubsection{Approximation by Low-Degree Trigonometric Polynomials}
A small Fourier ratio implies that the mollified measure $(f\mu)*\psi_{R^{-1}}$ has low analytic complexity.
\begin{itemize}
    \item \textbf{$L^2$ Approximation (Theorem \ref{theorem:L2approximation}):} For any $\eta > 0$, there exists a trigonometric polynomial $P$ with
    \[
    \|(f\mu)*\psi_{R^{-1}} - P\|_2 \leq \eta \|(f\mu)*\psi_{R^{-1}}\|_2
    \]
    and degree bounded by $\eta^{-2} (R^d \cdot \FR_{\mu,R}(f)^2 - 1)$. Similar results hold in $L^1$ (Theorem \ref{theorem:L1approximation}) and $L^\infty$ (Theorem \ref{Linfinityapproximation}).
    \item The approximation is constructive via a \textbf{random Fourier sampling scheme} that selects frequencies according to the distribution $|\widehat{(f\mu)_{R^{-1}}}|$.
\end{itemize}

\subsubsection{Restriction Theory and the Deterministic--Random Dichotomy}
The Fourier ratio sharply distinguishes between measures supporting classical restriction estimates and those arising from random fractal constructions.
\begin{itemize}
    \item \textbf{Deterministic Curvature (Circle):} For the arc-length measure $\sigma$ on $S^1$ and $f$ supported on an arc of length $R^{-1/2}$,
    \[
    \FR_{\sigma,R}(f) \lesssim R^{-1/4}
    \]
    (Theorem \ref{theorem:restrictionFR}), reflecting polynomial Fourier decay forced by curvature.
    \item \textbf{Random Fractal (Laba--Wang):} For their random Cantor measure $\mu$ of dimension $1$ in $\mathbb{R}^2$ and any $f \geq 0$,
    \[
    \FR_{\mu,R}(f) \gtrsim_{\epsilon} R^{-\epsilon} \quad \text{for every } \epsilon > 0
    \]
    (Theorem \ref{theorem:restrictionFR}), a consequence of its near-optimal extension estimates. This subpolynomial bound forces high approximation degree: $(f\mu)*\psi_{R^{-1}}$ \textbf{cannot} be approximated in $L^2$ with fixed accuracy by trigonometric polynomials of degree $o(R^2)$ (Proposition \ref{prop:labawanghighdegree}).
\end{itemize}

\subsubsection{Convex Geometry and the Complexity of Surface Measures}
For the surface measure $\mu$ on the boundary of a convex body $K \subset \mathbb{R}^d$, the Fourier ratio detects the geometric complexity of its normal set.
\begin{itemize}
    \item \textbf{Stationary Phase Analysis (Theorem \ref{theorem:convexstationaryphase}):} The Fourier transform $\widehat{\mu_{R^{-1}}}$ is essentially concentrated on directions near the set $N(K)$ of outward unit normals.
    \item \textbf{Approximation Degree via Minkowski Dimension (Corollary \ref{cor:minkowskikicksass}):} If $N(K)$ has upper Minkowski dimension $a$, then for any $\epsilon > 0$, the mollified measure $\mu_{R^{-1}}$ can be approximated in $L^2$ with relative error $\eta$ by a trigonometric polynomial of degree $\lesssim_{\eta,\epsilon} R^{a+1+\epsilon}$.
\end{itemize}

\subsection{Connections and Broader Context}
Our work builds a bridge between discrete Fourier analysis and the continuum theory of measures. The Fourier ratio $\FR_{\mu,R}(f)$ plays precisely the same role as its discrete counterpart: it controls uncertainty principles, measures pseudorandomness versus structure, and dictates the complexity of trigonometric approximation. The continuous setting naturally incorporates geometric information---covering numbers, Minkowski dimensions, and curvature---which enriches the discrete theory while preserving its essential philosophy.

The results also reveal a unified landscape connecting several active areas: fractal uncertainty principles \cite{BD18}, Fourier restriction theory \cite{Stein1993, LW18}, random trigonometric approximation \cite{Talagrand98}, and geometric measure theory on convex sets \cite{BCT97, BRT98}. The Fourier ratio serves as the common parameter through which these diverse phenomena interact.

\subsection{Structure of the Paper}
Section \ref{section:previousresults} reviews the discrete Fourier ratio theory. Section \ref{section:lowerandupperbounds} establishes the basic bounds and the fractal uncertainty principle. Section \ref{section:randomtrigapproximation} presents the trigonometric approximation theorems. Section \ref{section:fourierratioandrestriction} explores the deterministic--random dichotomy via restriction theory. Section \ref{section:convexgeometry} applies the theory to convex surfaces. Section \ref{section:discretecomparison} systematically compares the discrete and continuous theories. The results are summarized in Section \ref{section:summary}. All proofs are contained in Section \ref{section:proofs}.

\vskip.125in 

\subsection{Notation and Conventions}

We collect here the principal notation and conventions used throughout the paper.

\medskip

\noindent\textbf{Ambient space.}
We work in Euclidean space $\mathbb{R}^d$, $d \geq 1$. Lebesgue measure on $\mathbb{R}^d$
is denoted by $|\cdot|$.

\medskip

\noindent\textbf{Measures and functions.}
Throughout, $\mu$ denotes a compactly supported Borel measure on $\mathbb{R}^d$,
and $f \in L^2(\mu)$. The support of $f$ (with respect to $\mu$) is denoted by
\[
E_f := \operatorname{supp}(f).
\]

\medskip

\noindent\textbf{Fourier transform.}
For a finite Borel measure $\nu$ on $\mathbb{R}^d$, we use the normalization
\[
\widehat{\nu}(\xi) = \int e^{-2\pi i x \cdot \xi} \, d\nu(x).
\]
In particular, $\widehat{f\mu}$ denotes the Fourier transform of the measure
$f\,d\mu$.

\medskip

\noindent\textbf{Approximate identities.}
Let $\psi \in C_c^\infty(\mathbb{R}^d)$ be a fixed nonnegative function satisfying
$\int \psi = 1$. For $\delta > 0$, define
\[
\psi_\delta(x) := \delta^{-d}\psi(x/\delta).
\]
Throughout the paper we work at scale $\delta = R^{-1}$, where $R \geq 1$ is the
frequency parameter.

\medskip

\noindent\textbf{Regularized Fourier norms.}
For $1 \leq p \leq \infty$, we define the regularized Fourier $L^p$-norm at scale $R$ by
\[
X_{p,\mu,R}(f)
:=
\left(
R^{-d}
\int_{\mathbb{R}^d}
\big|
\widehat{f\mu} * \psi_{R^{-1}}(\xi)
\big|^p
\, d\xi
\right)^{1/p},
\]
with the usual modification when $p=\infty$.

\medskip

\noindent\textbf{Fourier ratio.}
The (continuous) Fourier ratio of $f$ with respect to $\mu$ at scale $R$ is
\[
\mathrm{\FR}_{\mu,R}(f)
:=
\frac{X_{1,\mu,R}(f)}{X_{2,\mu,R}(f)}.
\]

\medskip

\noindent\textbf{Neighborhoods and covering scale.}
For a set $E \subset \mathbb{R}^d$ and $r>0$, we denote by $E_r$ the $r$--neighborhood of
$E$. In particular,
\[
E_{R^{-1}}^f := (E_f)_{R^{-1}}.
\]
The quantity $R^d |E_{R^{-1}}^f|$ is comparable to the number of balls of radius
$R^{-1}$ needed to cover $E_f$.

\medskip

\noindent\textbf{Asymptotic notation.}
We write $A \lesssim B$ if $A \leq C B$ for a constant $C$ depending only on fixed
parameters such as the dimension. The notation $A \approx B$ means
$A \lesssim B \lesssim A$.

\section{Previous results} \label{section:previousresults} The discrete analog of the Fourier ratio above was used to study the complexity of time series in \cite{A2025}. We shall describe these in considerable detail, as the results have a direct impact on the investigations in this paper. The authors' starting point was the following result due to Bourgain and Talagrand (\cite{Talagrand98}). 

\begin{theorem}\label{thm:Talagrand98}
Let $(\varphi_j)_{j=1}^n$ be an orthonormal system in 
$L^2(\mathbb{Z}_N)$ with $\|\varphi_j\|_{L^\infty} \leq K$ for $1 \leq j \leq n$.  
There exists a constant $\gamma_0 \in (0,1)$ and a subset 
$I \subset \{1, \dots, n\}$ with $|I| \ge \gamma_0 n$  
such that for every $a = (a_i) \in \mathbb{C}^n$,
$$
\left( \sum_{i \in I} |a_i|^2 \right)^{1/2} 
\leq C_T \, K \, \big( \log(n) \, \log \log(n) \big)^{1/2} 
\left\| \sum_{i \in I} a_i \varphi_i \right\|_{L^1},
$$
where $C_T > 0$ is a universal constant.
\end{theorem}

We shall need a version of this, stated for signals on ${\mathbb Z}_N$, with the roles of $h$ and $\widehat{h}$ reversed (see \cite{BIMN2025}). First, we need a definition. 

\begin{definition} \label{def:generic} Let $0< p <1$. Then, a random set $S \subset [n]= \{0, 1, \cdots, n-1\}$ is generic if each element of $[n]$ is selected independently with probability $p$.
\end{definition} 

The following result can be deduced from Theorem \ref{thm:Talagrand98} (see \cite{BIMN2025}). 

\begin{theorem} \label{thm: combo} There exists $\gamma_0 \in (0,1)$ such that if $h: {\mathbb Z}_N \to {\mathbb C}$ supported in a generic set $M$ (in the sense of Definition \ref{def:generic}) of size $\gamma_0 \frac{N}{\log(N)}$, then with probability $1-o_N(1)$, 
\begin{equation} \label{eq: talagrand} {\left( \frac{1}{N}  \sum_{m \in {\mathbb Z}_N} {|\widehat{h}(m)|}^2 \right)}^{\frac{1}{2}} \leq C_T {(\log(N) \log \log(N))}^{\frac{1}{2}} \cdot \frac{1}{N} \sum_{m \in {\mathbb Z}_N} |\widehat{h}(m)|, \end{equation} where $C_T > 0$ is a constant that depends only on $\gamma_0$, and 
$$ \widehat{h}(m)=N^{-\frac{1}{2}} \sum_{x \in {\mathbb Z}_N} \chi(-xm) h(x), \ \chi(t) \equiv e^{\frac{2 \pi i t}{N}}.$$ 
\end{theorem} 

\vskip.125in 

\begin{remark} \label{rmk:talagrandsharpness} It is known (see \cite{Talagrand98}) that $\sqrt{\log(N)}$ in (\ref{eq: talagrand}) cannot, in general, be removed, and it is not known whether the removing the remaining $\sqrt{\log \log(N)}$ is possible. \end{remark} 

In some cases, the term ${(\log(N) \log \log(N))}^{\frac{1}{2}}$ can be removed. The following result, stated in this setting, is due to Bourgain \cite{Bourgain89}. 

\begin{theorem} \label{thm:bourgain} Suppose that $M$ is generic, as above, $|M|= \lceil N^{\frac{2}{q}} \rceil$, $q>2$. Then for all $h: {\mathbb Z}_N \to {\mathbb C}$, supported in $M$, 
\begin{equation} \label{eq:bourgain} {\left( \frac{1}{N} \sum_{m \in {\mathbb Z}_N} {|\widehat{h}(m)|}^q \right)}^{\frac{1}{q}} \leq C(q) \cdot {\left( \frac{1}{N} \sum_{m \in {\mathbb Z}_N} {|\widehat{h}(m)|}^2 \right)}^{\frac{1}{2}} \end{equation} 

It is not difficult to deduce from the proof that if $M$ is chosen randomly, then the bound holds with probability $1-\epsilon$ if $C(q)$ is replaced by $\frac{C(q)}{\epsilon}$. 
\end{theorem} 

\vskip.125in 

We learned the following observation from William Hagerstrom \cite{Hagerstrom2025}, which can be proven using H\"older's inequality. 

\begin{lemma} \label{lemma:vershynintrick} Suppose that for $h: {\mathbb Z}_N \to {\mathbb C}$, 
\begin{equation} \label{eq:bourgainproxy} {\left( \frac{1}{N} \sum_{x \in {\mathbb Z}_N} {|\widehat{h}(x)|}^q \right)}^{\frac{1}{q}} \leq C(q) {\left( \frac{1}{N} \sum_{x \in {\mathbb Z}_N} {|\widehat{h}(x)|}^2 \right)}^{\frac{1}{2}} \end{equation} for some $q>2$. 

Then 
$$ {\left( \frac{1}{N} \sum_{x \in {\mathbb Z}_N} {|\widehat{h}(x)|}^2 \right)}^{\frac{1}{2}} \leq {(C(q))}^{\frac{q}{q-2}} \cdot \frac{1}{N} \sum_{x \in {\mathbb Z}_N} |\widehat{h}(x)|. $$
\end{lemma} 

\begin{remark} \label{rmk:nolog} It follows that if the expected size of $M$ in Theorem \ref{thm: combo} is $O(N^{1-\epsilon})$ for some $\epsilon>0$, then 
$C_T \sqrt{\log(N) \log \log(N)}$ in (\ref{eq:bourgain}) can be replaced by $C'_T$, independent of $N$. \end{remark} 

\vskip.125in 

Let's explore the setup above a bit. Let $S \subset {\mathbb Z}_N$, chosen randomly. Then in view of Theorem \ref{thm: combo} and Remark \ref{rmk:nolog}, with very high probability, we have 
\begin{equation} \label{eq:mamaindicator} {\left( \frac{1}{N} \sum_{x \in {\mathbb Z}_N} {|\widehat{1}_S(x)|}^2 \right)}^{\frac{1}{2}} \leq C_T \frac{1}{N} \sum_{x \in {\mathbb Z}_N}|\widehat{1}_S(x)|.\end{equation}  

Using Plancherel, we see that 
\begin{equation}\label{eq:mamaindicatorconcretelowerbdd}
\frac{1}{C_T} \sqrt{\frac{|S|}{N}} \leq \frac{1}{N} \sum_{x \in {\mathbb Z}_N}|\widehat{1}_S(x)|.
\end{equation}

On the other hand, let's consider the case where $S$ is far from random. Suppose that $S=\{0\}$. Then $\widehat{1}_S(x) \equiv \frac{1}{\sqrt{N}}$ and it is not difficult to check that in this case, we have $\frac{1}{N} \sum_{x \in {\mathbb Z}_N}|\widehat{1}_S(x)| = \frac{1}{\sqrt{N}} \cdot \sqrt{\frac{1}{N}},$ so for a fixed $S\subseteq \mathbb{Z}_N$, $\frac{1}{N} \sum_{x \in {\mathbb Z}_N}|\widehat{1}_S(x)| = \frac{1}{\sqrt{N}} \cdot \sqrt{\frac{|S|}{N}},$ a much smaller quantity than what we get in (\ref{eq:mamaindicatorconcretelowerbdd}) in the generic case. A similar calculation can be carried out in the case when $f$ is the indicator function of a subgroup of ${\mathbb Z}_N$. This suggests that the quantity $\frac{\frac{1}{N} \sum_{x \in {\mathbb Z}_N}|\widehat{1}_S(x)|}{{\left( \frac{1}{N} \sum_{x \in {\mathbb Z}_N} {|\widehat{1}_S(x)|}^2 \right)}^{\frac{1}{2}}}$ may contain information indicating the degree to which a set $S$ is random. 

More generally, under the assumptions of Theorem \ref{thm:bourgain}, with the additional assumption that $|M|=O(N^{1-\delta})$ for some $\delta>0$, we have 
\begin{equation} \label{eq:randomfunction} \frac{1}{C_T} \leq \frac{\frac{1}{N} \sum_{x \in {\mathbb Z}_N}|\widehat{f}(x)|}{{\left( \frac{1}{N} \sum_{x \in {\mathbb Z}_N} {|\widehat{f}(x)|}^2 \right)}^{\frac{1}{2}}}, \end{equation} and, once again, the examples above suggest that it may be reasonable to use this quantity as an indicator of the degree of randomness of the signal $f$. Since the constant $C_T$ is of practical significance, we conduct some numerical experiments below to understand what its value may be. 

\vskip.125in 

\subsubsection{Bounding the Fourier ratio from below.}

Suppose that 
$$\frac{1}{N} \sum_{x \in {\mathbb Z}_N} |\widehat{f}(m)| \leq \epsilon {\left( \frac{1}{N} \sum_{x \in {\mathbb Z}_N}{|\widehat{f}(m)|}^2 \right)}^{\frac{1}{2}}.$$ By the triangle inequality, 
$$ |f(x)| \leq N^{-\frac{1}{2}} \sum_{m \in \mathbb{Z}_N} |\widehat{f}(m)|=N^{\frac{1}{2}} \cdot \frac{1}{N} \sum_{m \in {\mathbb Z}_N} |\widehat{f}(m)| .$$ 
$$ \leq \sqrt{N} \cdot \epsilon \cdot {\left( \frac{1}{N} \sum_{m \in {\mathbb Z}_N}{|\widehat{f}(m)|}^2 \right)}^{\frac{1}{2}}.$$

Squaring both sides, summing over ${\mathbb Z}_N$ and taking square roots yields $\epsilon \ge \frac{1}{\sqrt{N}}$, and one can see this bound is realized by the constant function $1$. We just established that 
$$ \frac{1}{\sqrt{N}} \leq \frac{\frac{1}{N} \sum_{x \in {\mathbb Z}_N}|\widehat{f}(x)|}{{\left( \frac{1}{N} \sum_{x \in {\mathbb Z}_N} {|\widehat{f}(x)|}^2 \right)}^{\frac{1}{2}}} \leq 1,$$ where the upper bound follows by Cauchy-Schwarz. 

\vskip.125in 

The same argument shows that if $f$ is supported in $E \subset {\mathbb Z}_N$, then 
\begin{equation} \label{eq:discreteratiobounds} \frac{1}{\sqrt{|E|}} \leq \frac{\frac{1}{N} \sum_{x \in {\mathbb Z}_N}|\widehat{f}(x)|}{{\left( \frac{1}{N} \sum_{x \in {\mathbb Z}_N} {|\widehat{f}(x)|}^2 \right)}^{\frac{1}{2}}}. \end{equation} 

To see that (\ref{eq:discreteratiobounds}) can be realized, let $E={\mathbb Z}_p$, $p$ prime, sitting inside ${\mathbb Z}_{pq}$ in the obvious way. Then, setting $N=pq$, we have 
$$ \widehat{1}_E(m)=\frac{1}{\sqrt{N}} \sum_{k=0}^{p-1} e^{-\frac{2 \pi i km}{p}} =\frac{p}{\sqrt{N}} 1_S(m),$$ where $S$ is ${\mathbb Z}_q$, sitting inside ${\mathbb Z}_{pq}$ in the natural way. Setting $f=1_E$, it follows that 
$$  \frac{\frac{1}{N} \sum_{x \in {\mathbb Z}_N}|\widehat{f}(x)|}{{\left( \frac{1}{N} \sum_{x \in {\mathbb Z}_N} {|\widehat{f}(x)|}^2 \right)}^{\frac{1}{2}}}=\frac{1}{\sqrt{p}}=\frac{1}{\sqrt{|E|}}.$$

\vskip.125in 

\subsubsection{The main results of \cite{A2025}} The first result shows that if $f: {\mathbb Z}_N \to {\mathbb C}$ is concentrated on a random set, then the Fourier ratio is very large with high probability. 

\begin{theorem} \label{thm:concentration} Let $f: {\mathbb Z}_N \to {\mathbb C}$. Suppose that there exists a generic set $M$ such that
$$ {\|f\|}_{L^2(M^c)} \leq r {\|f\|}_2$$ for some $r \in (0,1)$, with $|M| \leq \gamma_0 \frac{N}{\log(N)}$, where $\gamma_0$ is as in Theorem \ref{thm: combo}. Suppose that 
\begin{equation} \label{eq:concentrationrestriction} r< \frac{1-r}{C_T \sqrt{\log(N) \log \log(N)}}. \end{equation}

Then 
\begin{equation} \label{eq:largeconcentrationratio}  \frac{1-r\frac{C_T \sqrt{\log(N) \log \log(N)}}{1-r}}{\frac{C_T \sqrt{\log(N) \log \log(N)}}{1-r}} \leq \frac{\frac{1}{N} \sum_{x \in {\mathbb Z}_N}|\widehat{f}(x)|}{{\left( \frac{1}{N} \sum_{x \in {\mathbb Z}_N} {|\widehat{f}(x)|}^2 \right)}^{\frac{1}{2}}}\end{equation} with probability $1-o_N(1)$. 
\end{theorem} 

\vskip.125in 

Using the observation in Remark \ref{rmk:nolog}, we can replace (\ref{eq:concentrationrestriction}) with $r<\frac{1-r}{C_T}$, and we can replace (\ref{eq:largeconcentrationratio}) with 
$$ \frac{\left(1-r\frac{C_T }{1-r} \right)}{\frac{C_T }{1-r}} \leq \frac{\frac{1}{N} \sum_{x \in {\mathbb Z}_N}|\widehat{f}(x)|}{{\left( \frac{1}{N} \sum_{x \in {\mathbb Z}_N} {|\widehat{f}(x)|}^2 \right)}^{\frac{1}{2}}}.$$ 

\vskip.125in 

The next result established in \cite{A2025} shows that if the Fourier ratio is suitably small, then the signal can be well-approximated by a trigonometric polynomial of a low degree. 

\begin{theorem}\label{thm:L2polynomialapprox} Let $f:\mathbb Z_N \to \mathbb C$, and let $\eta>0$. Let 
$$ \epsilon = \frac{\frac{1}{N} \sum_{x \in {\mathbb Z}_N}|\widehat{f}(x)|}{{\left( \frac{1}{N} \sum_{x \in {\mathbb Z}_N} {|\widehat{f}(x)|}^2 \right)}^{\frac{1}{2}}}.$$

Then for any $k$ such that $$k > \frac{N\epsilon^2 - 1}{\eta^2},$$ there exists a trigonometric polynomial
$$P(x) = \sum_{i=1}^k c_i \chi(m_ix)$$ such that
$$\|f - P\|_2 < \eta \|f\|_2.$$
\end{theorem}

\vskip.125in 

The next result establishes this type of an approximation in the $L^{\infty}$ norm. 

\begin{theorem}\label{thm:Linftypolynomialapprox} Let $f:\mathbb Z_N\to\mathbb C$ and let $\eta>0$. Then for any $k$ such that $$k > 8\left(\frac{\|\widehat f\|_{L^1(\mu)}}{\|f\|_\infty}\right)^2 \frac{N\log (4N)}{\eta^2},$$
there exists a trigonometric polynomial
$$P(x) = \sum_{i=1}^k c_i \chi(m_ix)$$ such that
$$\|f - P\|_\infty < \eta \|f\|_\infty.$$
\end{theorem}

\vskip.125in 

\begin{remark} Note that the triangle inequality shows $\frac{\|\widehat f\|_{L^1(\mu)}}{\|f\|_\infty} \geq N^{-\frac12}$, and so in the best case, Theorem \ref{thm:Linftypolynomialapprox} indeed gives a polynomial of degree $O(\log(N))$.\end{remark}

\vskip.125in 

The following result demonstrates that the discrete Fourier ratio serves as the controlling parameter in the classical uncertainty principle. 

\begin{theorem} \label{theorem:FRdiscreteuncertainty} Let $f: {\mathbb Z}_N \to {\mathbb C}$, $L^2$-concentrated in $E \subset {\mathbb Z}_N$ at level $a \in (0,1)$, in the sense that 
$$ {||f||}_{L^2(E^c)} \leq a {||f||}_{L^2({\mathbb Z}_N)},$$ with $\widehat{f}$
$L^1$-
concentrated
on $S \subset {\mathbb Z}_N$ at level $b\in (0,1)$, in the sense that 
$$ {||\widehat{f}||}_{L^1(S^c)} \leq b{||\widehat{f}||}_{L^1({\mathbb Z}_N)}.$$ 

Then 
\begin{equation}\label{eq:FRcontrollingparameter}  (1-a)^2  \cdot \frac{N}{|E|} \leq {\FR(f)}^2 \leq \frac{|S|}{{(1-b)}^2}. \end{equation} 

\end{theorem} 

\vskip.125in 

In particular, 
\begin{equation} \label{eq:FRintermediary} {(1-a)}^2 \cdot {(1-b)}^2 \cdot N \leq |E| \cdot |S|, \end{equation} a version of the classical uncertainty principle (see e.g. \cite{DS89}). While (\ref{eq:FRintermediary}) is well-known, (\ref{eq:FRcontrollingparameter}) shows that the Fourier Ratio is a natural controlling parameter in the Fourier Uncertainty Principle.  

The main thrust of this paper is to understand the degree to which these ideas generalize to the rather more complicated setting of compactly supported measures in Euclidean space. 

\vskip.125in 

\section{Lower and upper bounds on the Fourier ratio} \label{section:lowerandupperbounds}

Our first result is a continuous analog of (\ref{eq:discreteratiobounds}), a lower bound on $\FR(f\mu)$. 

\begin{proposition} \label{prop:lowerbound} Let $\mu$ be a compactly supported Borel measure, and let $f \in L^1 \cap L^2(\mu)$. Let $E_f$ denote the support of $f$. Then 
\begin{equation} \label{eq:continuousratiobounds} \sqrt{\frac{1}{R^d|E_f^{\frac{1}{R}}|}} \leq \frac{X_{1,\mu, R}(f)}{X_{2,\mu, R}(f)}, \end{equation} where $E_f^{\frac{1}{R}}$ denotes the $\frac{1}{R}$-neighborhood of the support of $f$. 
\end{proposition} 

\vskip.125in 

\begin{remark} One way to see how Proposition \ref{prop:lowerbound} is analogous to (\ref{eq:discreteratiobounds}) is by noting that 
$$ |E^{\frac{1}{R}}| \sim R^{-d} \cdot \# \ \text{balls of radius} \ R^{-1} \ \text{needed to cover} \ E.$$ It follows that the right hand side of (\ref{eq:continuousratiobounds}) can be replaced by 
$$ \sqrt{\frac{1}{\# \ \text{balls of radius} \ R^{-1} \ \text{needed to cover} \ E}},$$ analogously to (\ref{eq:discreteratiobounds}). 
\end{remark} 
 
\vskip.125in 

If $\mu$ is the surface measure on a convex hypersurface in ${\mathbb R}^d$, $d \ge 2$, Proposition \ref{prop:lowerbound} shows that 
$$ \frac{c}{R^{\frac{d-1}{2}}} \leq \frac{X_{1,\mu,R}(1)}{X_{2,\mu,R}(1)}, $$ which shows that polyhedra nearly optimize this inequality, in view of (\ref{eq:convexratio}). 

\vskip.125in 

\subsection{An upper bound on the Fourier ratio and a fractal uncertainty principle} \label{subsection:uncertaintyprinciple} We now get an upper bound that complements the result we have obtained above. To this end, we have the following result. 

\begin{theorem} \label{theorem:fractalsandwich} Suppose that $\widehat{{(f\mu)}_{R^{-1}}}$ is $L^1$-concentrated in $X \subset {\mathbb R}^d$ in the sense that for some $\eta \in (0,1)$, 
\begin{equation} \label{eq:L1concentration} {||\widehat{f\mu_{R^{-1}}}||}_{L^1(X_R^c)} \leq \eta \cdot {||\widehat{f\mu_{R^{-1}}}||}_{L^1},\end{equation} where $X_R=X \cap B_{100R}$. Then 
\begin{equation} \label{eq:FRsandwich} \sqrt{\frac{1}{R^d|E_f^{\frac{1}{R}}|}} \leq \FR(f\mu) \leq 
\sqrt{\frac{|X_R|}{R^d {(1-\eta)}^2}}. \end{equation} 

It follows that 
\begin{equation} \label{eq:fractaluncertainty} {(1-\eta)}^2 \leq |E_f^{R^{-1}}| \cdot |X_R|. \end{equation}  

\vskip.125in 

If we assume that there exist $c, C_X$ universal constants, and $s_f, \alpha_X \in (0,d)$, such that 
\begin{equation} \label{eq:porosity} |X_R| \leq C_X R^{\alpha_X}, \end{equation} and 
\begin{equation} \label{eq:sizespacesupport} |E_f^{R^{-1}}| \leq cR^{-d+s_{f}}, \end{equation} for all large enough $R$, then in order to ensure the condition (\ref{eq:fractaluncertainty}) for all large enough $R$, we can ask for the stronger condition 
$$ d < s_{f}+\alpha_X.$$
\end{theorem} 

\vskip.125in 

\subsubsection{Geometric interpretation of \eqref{eq:sizespacesupport} and \eqref{eq:porosity}}

We will interpret equation \eqref{eq:sizespacesupport} geometrically:

The upper box dimension of $E\subseteq \R^d$ can be defined as \[\overline{\dim}_B (E):=d-\liminf_{R\to \infty} \frac{\log (|E^{R^{-1}}|)}{\log (R^{-1})}.\]

 Observe that, if $|E_f^{R^{-1}}|\leq c R^{s_f-d}$ for all $R\geq R_0$, then $\overline{\dim}_B (E)\leq s_f$. On the other hand, if $\overline{\dim}_B (E)< s_f$, then there is a large $R_0$ so that $|E_f^{R^{-1}}|\leq R^{s_f-d}$ for all $R\geq R_0$.

Hence, equation \eqref{eq:sizespacesupport} can be replaced by the slightly stronger condition $\overline{\dim}_B (E)< s_f$.

We can also interpret equation \eqref{eq:porosity} geometrically:

We define the \textit{logarithmic upper asymptotic density of} $X$ as \[\overline{d}_{\log}(X):=\limsup_{R\to \infty} \frac{\log(X\cap B_{100R})}{\log(R)}.\]

Note that the definition is unchanged if we replace the numerator by $\log(X \cap B_R)$ which is the expression one would naturally use to define it.

Observe that if $|X_R|:=|X \cap B_{100R}|\leq C_X R^{\alpha_X}$, then $\overline{d}_{\log}(X)\leq \alpha_X$. On the other hand, if  $\overline{d}_{\log}(X)< \alpha_X$, then $|X_R|:=|X \cap B_{100R}|< R^{\alpha_X}$ for all large enough $R$.

Hence, equation \eqref{eq:porosity} can be replaced by the slightly stronger condition $\overline{d}_{\log}(X)< \alpha_X$.

\vskip.125in 

\subsection{$L^1$ concentration cannot be replaced by $L^2$ concentration}

In Theorem~\ref{theorem:fractalsandwich}, the upper bound on the Fourier ratio
$$
\FR(f\mu)=\frac{X_{1,\mu,R}(f)}{X_{2,\mu,R}(f)}
$$
depends essentially on the $L^1$ concentration assumption
$$
\|\widehat{(f\mu)_{R^{-1}}}\|_{L^1(S_R^c)}\le \eta \|\widehat{(f\mu)_{R^{-1}}}\|_{L^1}.
$$
The purpose of this subsection is to show that this hypothesis cannot be replaced by $L^2$ concentration.  
In fact, even if $\widehat{(f\mu)_{R^{-1}}}$ is almost completely supported in $S_R$ in the $L^2$ sense, no nontrivial upper bound on $\FR(f\mu)$ follows, and the uncertainty inequality of Theorem~\ref{theorem:fractalsandwich} fails completely. We have the following continuous variant of the discrete variant in Theorem \ref{theorem:FRdiscreteuncertainty}. 

\begin{proposition} \label{prop:L2fails}
Let $\mu$ be Lebesgue measure on $\mathbb R^d$, and let $R \ge 1$. Let $S_R \subset \mathbb R^d$ be an open ball of radius $R$ and finite positive measure.

Then for every constant $C>0$ there exist $b \in (0,1)$ and a function $f \in L^1(\mu) \cap L^2(\mu)$ such that
$$
\|\widehat{(f\mu)_{R^{-1}}}\|_{L^2(S_R^c)} \le b \|\widehat{(f\mu)_{R^{-1}}}\|_{L^2(\mathbb R^d)},
$$
but
$$
\FR(f\mu) > \frac{C}{1-b} |S_R|^{1/2} R^{-d/2}.
$$

In particular, even for this fixed measure $\mu$, no estimate of the form
$$
\FR(f\mu) \le \frac{C}{1-b} |S_R|^{1/2} R^{-d/2}
$$
can be deduced from the $L^2$ concentration condition alone, with a constant $C$ independent of $f$ and $b$. Thus the $L^1$ concentration hypothesis in Theorem \ref{theorem:fractalsandwich} cannot, in general, be weakened to $L^2$ concentration.
\end{proposition}

\vskip.125in 

\subsection{Applications to signal recovery} \label{section:signalrecovery} In this subsection, we are going to apply the uncertainty principle derived in Theorem \ref{theorem:fractalsandwich} to the problem of signal recovery in the setting of compactly supported Borel measures on ${\mathbb R}^d$. 

\begin{theorem} \label{theorem:measurerecovery} Let $\mu$ be a compactly supported Borel measure on a set $E \subset {\mathbb R}^d$ such that 
\begin{equation} \label{eq:spacesize} |E^{R^{-1}}| \leq C_E R^{s_E -d} \end{equation} for some $C_E$ independent of $R>1$, where $E^{R^{-1}}$ is the $R^{-1}$ neighborhood of $E$. 

Suppose that $X \subset {\mathbb R}^d$ such that 
\begin{equation} \label{eq:spectrumsize} |X \cap B_R| \leq C_XR^{\alpha_X}, \end{equation} for some $\alpha_X\geq 0$ and $C_X$ independent of $R>1$, where $B_R$ is the ball of radius $R$ centered at the origin. 

Suppose that the frequencies ${\{\widehat{\mu}(\xi)\}}_{\xi \in X \cap B_R}$ are unobserved, where $B_R$ is the ball of radius $R>0$ centered at the origin. Then $\mu$ can be recovered exactly and uniquely provided that 
$$ s_E+\alpha_X <d$$ and R is sufficiently large.
\end{theorem} 

\vskip.125in 

\begin{remark} We prove Theorem \ref{theorem:measurerecovery} using Logan's celebrated $L^1$-minimization idea. The interested reader can check that we could have proven it using the uncertainty principle (\ref{eq:fractaluncertainty}) that follows the classical argument due to Donoho and Stark in the discrete setting (\cite{DS89}). We leave the details to the interested reader. 
\end{remark} 

\vskip.25in 

\section{Approximation via random Fourier sampling} \label{section:randomtrigapproximation} 

We are now going to see that a small Fourier ratio implies that $f \in L^2(\mu)$ can be well approximated by a trigonometric polynomials of a low degree. 

\begin{theorem} \label{theorem:L2approximation} Let $\mu$ denote the restriction of the $s$-dimensional Hausdorff measure to $E$. Let $f \in L^2(\mu_E)$. Then there exists a trigonometric polynomial $P$ such that 
\begin{equation} \label{eq:L2approximation} {||(f\mu)*\psi_{R^{-1}}-P||}_2 \leq \eta \cdot {||(f\mu)*\psi_{R^{-1}}||}_2. \end{equation} 

Moreover, the degree of this polynomial is 
\begin{equation} \label{eq:degree} \leq \frac{1}{\eta^2} \left( \frac{R^d X^2_{1,\mu,R}(f)}{X^2_{2,\mu,R}(f)} -1 \right). \end{equation} 


\end{theorem} 

\vskip.125in 

\begin{remark} (Polynomial degree count in Theorem \ref{theorem:L2approximation})\label{remark:L2interpretation} In view of Proposition \ref{prop:lowerbound}, the degree of the approximating polynomial in (\ref{eq:degree}) is at most 
\begin{equation} \label{eq:highestL2degree} C\frac{R^d}{\eta^2}\end{equation} and at least 
$$ \frac{1}{\eta^2} \cdot \frac{1}{|E_f^{R^{-1}}|},$$ where, as before, $E_f$ is the support of $f$, $E_f^{R^{-1}}$ denotes its $R^{-1}$ neighborhood, and $|E_f^{R^{-1}}|$ is the Lebesgue measure of this neighborhood. 

The lower bound shows that the smallest possible degree of the polynomial approximating $(f\mu)*\psi_{R^{-1}}$ is 
\begin{equation} \label{eq:lowestL2degree} \approx \frac{1}{\eta^2} \cdot \frac{R^d}{\# \ \text{balls of radius} \ R^{-1} \ \text{needed to cover the} \ support(f)}. \end{equation} 

It is not difficult to check that (\ref{eq:lowestL2degree}) is realized if $\mu$ is the natural measure on a compact piece of a $k$-plane in $\mathbb{R}^d$, $1 \leq k \leq d-1$. The highest possible degree, given by (\ref{eq:highestL2degree}) is realized, for example, by the natural measure on the unit sphere in ${\mathbb R}^d$, $d \ge 2$, with $f \equiv 1$. We give a quick sketch. Let $\sigma$ denote the surface measure on $S^{d-1}$. A direct calculation using the formula for the Fourier transform of the surface measure on the sphere (see e.g. \cite{Stein1993}, \cite{Strichartz83}) shows that 
$$ \FR(f\mu) \ge c>0,$$ so (\ref{eq:highestL2degree}) yields $\sim \frac{R^d}{\eta^2}$ as the upper bound on the degree of the approximating polynomial. It is not difficult to see that this degree be substantially lowered because the Fourier transform of $\sigma$ restricted to a $R^{-1}$ ball and convolved with the approximation to the identity at level $R^{-1}$ is concentrated on a cone of length $\approx R$ and aperture $\frac{1}{R}$ pointing in the direction normal to the sphere at center of the ball. 
\end{remark} 

\vskip.125in 

\begin{remark} \label{remark:polygon} We have pointed out above that if $\mu$ is the arc-length measure, say, on a convex polygon with finitely many sides, then the Fourier ratio is $\approx R^{-\frac{1}{2}} \log(R)$, where the implicit constant depends on the number of sides. Applying Theorem \ref{theorem:L2approximation}, this leads to the approximating polynomial of degree $\approx \log^2(R)$. At first glance, this may seem counterintuitive because one should be able to approximate the function $1$ on a convex polygon with finitely many sides with a trigonometric polynomial with finitely many terms. However, it is important to keep in mind that the input is the Fourier ratio, not the polygon itself. It is not difficult to see that if we consider a polygon inscribed in the unit circle, where the sides accumulate in a lacunary fashion at one point, then the Fourier ratio is, once again, $\approx R^{-\frac{1}{2}} \log^2(R)$. This shows that a polynomial of logarithmic degree is reasonable in this setting. 
\end{remark} 

\vskip.125in 

The next result establishes an $L^{\infty}$ version of Theorem \ref{theorem:L2approximation}. 

\begin{theorem}\label{Linfinityapproximation}

Let $\mu$ denote the restriction of the $s$-dimensional Hausdorff measure to $E$. Let $f\in L^{\infty}(\mu_E)$. Then there exists a trigonometric polynomial $P$ such that
\begin{equation}
\big\|(f\mu)*\psi_{R^{-1}} - P\big\|_{\infty} \le \eta\big\|(f\mu)*\psi_{R^{-1}}\big\|_{\infty}.
\end{equation}
Moreover, the degree of this polynomial can be chosen to satisfy
\begin{equation}
    \leq \dfrac{{32||\widehat{{(f\mu)*\psi_{R^{-1}}}||}_1^2}}{\eta^2\|(f\mu)*\psi_{R^{-1}}\|_\infty^2}\left(\log4C_d+d\log\left(\dfrac{8\pi R\|\widehat{{(f\mu)*\psi_{R^{-1}}}}\|_1}{\eta\|\widehat{{(f\mu)*\psi_{R^{-1}}}}\|_\infty}\right)\right)
\end{equation}



\end{theorem} 

\vskip.125in 

\begin{remark}(polynomial degree count in Theorem \ref{Linfinityapproximation})\label{remark:Linfinityinterpretation} Cauchy-Schwarz inequality and Plancherel inequality shows:
$$\frac{\|\widehat{(f\mu)*\psi_{R^{-1}}}\|_1}{\|f\|_\infty} \leq R^\frac{d}{2}|E_{f}^{R^{-1}}|^\frac{1}{2},$$ 
So if we plug in, the degree of the approximating polynomial has an upper bound 
$$\approx \frac{R^d|E_{f}^{R^{-1}}|}{\eta^2}\log\left(\frac{R^d|E_{f}^{R^{-1}}|}{\eta^\frac{2d}{d+2}}\right).$$ 

In our case, omitting the constant that only depends on dimension d, we see that the degree is at most $$R^d|E_{f}^{R^{-1}}| \approx \# \ \text{ balls of radius} \ R^{-1} \ \text{needed to cover the support of} \ f,$$ up to a logarithmic factor. 
\end{remark}

\vskip.125in 

\begin{theorem} \label{theorem:L1approximation} Let $\mu$ denote the restriction of the $s$-dimensional Hausdorff measure to $E$. Let $f \in L^2(\mu_E)$. Then there exists a trigonometric polynomial $P$ such that 
\begin{equation} {||(f\mu)*\psi_{R^{-1}}-P||}_1 \leq \eta \cdot {||(f\mu)*\psi_{R^{-1}}||}_1. 
\end{equation} 

Moreover, the degree of this polynomial is 
\begin{equation} \leq \frac{1}{\eta^2} \left( \frac{||\widehat{(f\mu)*\psi_{R^{-1}}}||_1}{{||(f\mu)*\psi_{R^{-1}}||}_1} \right)^2. 
\end{equation}



\end{theorem}

\vskip.125in 

\subsection{An exposition of the main ideas behind the proofs of the results in this section} 

We briefly describe the main idea behind the trigonometric polynomial approximation arguments in this section. For convenience, write 
$$ g(x)={(f\mu)*\psi_{R^{-1}}}(x), $$
so that
$$ g(x)=\int \widehat{g}(\xi) e^{2\pi i x \cdot \xi} d\xi. $$

The first step is to regard this Fourier inversion formula as an expectation of a random exponential. We define a probability measure on frequency space by 
$$ d{\mathbb P}(\xi)=\frac{|\widehat{g}(\xi)|}{{||\widehat{g}||}_1} d\xi, $$
and for each $\xi$ we associate the random variable
$$ Z(x)={||\widehat{g}||}_1 \operatorname{sgn}(\widehat{g}(\xi)) e^{2\pi i x \cdot \xi}. $$
A direct computation shows that for every $x$,
$$ {\mathbb E}(Z(x))=g(x). $$
In other words, $g$ is the expectation of a random character $Z$, with frequency distributed according to $|\widehat{g}|$.

We then take $k$ independent copies $Z_1,\dots,Z_k$ of $Z$ and set
$$ P(x)=\frac{1}{k} \sum_{j=1}^k Z_j(x). $$
The function $P$ is a random trigonometric polynomial of degree at most $R$ (up to constants depending on $\psi$), and by construction 
$$ {\mathbb E}(P(x))=g(x) \quad \text{for all } x. $$
The error $g-P$ is controlled by the variance of $Z$. For the $L^2$ result in Theorem \ref{theorem:L2approximation}, one computes
$$ {\mathbb E}(|g(x)-P(x)|^2)=\frac{1}{k} \operatorname{Var}(Z(x)) $$
and integrates over $x$. This yields
$$ {\mathbb E}{||g-P||}_2^2=\frac{1}{k}\left({||\widehat{g}||}_1^2 \cdot |E_f^{R^{-1}}| - {||g||}_2^2\right), $$
or the analogous expression over the whole space. Expressing ${||\widehat{g}||}_1$ and ${||\widehat{g}||}_2$ in terms of $X_{1,\mu,R}(f)$ and $X_{2,\mu,R}(f)$ shows that
$$ {\mathbb E}{||g-P||}_2^2 \le \frac{1}{k}\left(R^d X_{1,\mu,R}(f)^2 - X_{2,\mu,R}(f)^2\right), $$
which leads to the degree bound in Theorem \ref{theorem:L2approximation} after choosing 
$$ k \approx \eta^{-2}\left( \frac{R^d X_{1,\mu,R}(f)^2}{X_{2,\mu,R}(f)^2} -1\right). $$

For the $L^\infty$ result in Theorem \ref{Linfinityapproximation}, we use the same random construction, but we apply Hoeffding's inequality (or a similar concentration inequality) at each fixed spatial point $x$ to control $|P(x)-g(x)|$ with high probability. Since $\widehat{g}$ is supported in $\{|\xi|\lesssim R\}$, both $g$ and $P$ are Lipschitz with Lipschitz constant bounded by a constant times $R {||\widehat{g}||}_1$. We cover the relevant region in space by a $\delta$-net $N_\delta$ with 
$$ |N_\delta| \le C_d \frac{|E_f^{R^{-1}}|}{\delta^d}, $$
and combine pointwise concentration on the net with the Lipschitz property and a union bound. Choosing $\delta$ proportional to $\eta/(R {||\widehat{g}||}_1)$ leads to the logarithmic factors in the degree bound, and eventually yields the $L^\infty$ approximation with degree estimated in terms of ${||\widehat{g}||}_1$, ${||g||}_\infty$ and $|E_f^{R^{-1}}|$ as in Theorem \ref{Linfinityapproximation}.

The $L^1$ approximation in Theorem \ref{theorem:L1approximation} again starts from the same random polynomial $P$. We use the $L^2$ variance estimate together with Jensen's inequality and Hölder's inequality in the spatial variable to deduce that
$$ {\mathbb E}{||g-P||}_1 \le \frac{|E_f^{R^{-1}}|^{1/2}}{\sqrt{k}} {||\widehat{g}||}_1, $$
and then choose $k$ large enough so that this is at most $\eta {||g||}_1$. As before, the degree bound is naturally expressed in terms of the ratio ${||\widehat{g}||}_1/{||g||}_1$, which is a Fourier ratio type quantity.

In all three cases, the key point is that the approximation is obtained by sampling from the Fourier transform of $(f\mu)*\psi_{R^{-1}}$, and the number of frequencies required is governed by the relationship between the $L^1$ and $L^2$ sizes of the Fourier transform, encapsulated by the Fourier ratio $\FR(f\mu)$.

\vskip.25in 

\section{Restriction theory and the Fourier ratio} \label{section:fourierratioandrestriction}

In this section, we are going to explore the relationship between the Fourier Ratio $\FR(f\mu)$ and the restriction phenomenon. The main focus is to expose the difference between discrete and deterministic cases. 

Let us first remind ourselves of the Stein-Tomas restriction theorem for the circle. See \cite{Tomas1975}, and, for the endpoint case, see e.g. \cite{Stein1993}. 

\begin{theorem} \label{theorem:steintomas} Let $\sigma$ denote the arc-length measure on the circle. Then 
\begin{equation} \label{eq:steintomas} {||\widehat{g\sigma}||}_{L^6({\mathbb R}^2)} \leq C {||g||}_{L^2(\sigma)}. \end{equation} 
\end{theorem} 

A beautiful example of a restriction theorem for a random set in this setting is the following far-reaching result due to Izabella Laba and Hong Wang (\cite{LW18}). 

\begin{theorem} \label{theorem:labawang} 
Let $d \in \mathbb{N}$ and $0 < \alpha < d$. Then there exists a probability measure supported on a subset of $[0,1]^d$ of Hausdorff dimension $\alpha$ such that:
\begin{itemize}
    \item i) For every $0 < \gamma < \alpha$, there is a constant $C_1(\gamma)$ such that
    $$
        \mu(B(x,r)) \leq C_1(\gamma) r^\gamma 
        \quad \forall x \in \mathbb{R}^d, \; r > 0.$$ 
    
    \item ii) For every $\beta < \min(\alpha/2, 1)$, there is a constant $C_2(\beta) > 0$ such that
    $$ 
        |\widehat{\mu}(\xi)| \leq C_2(\beta)(1+|\xi|)^{-\beta}
        \quad \forall \xi \in \mathbb{R}^d.      
    $$

    \item iii) For every $\epsilon>0$ there is $C_{\epsilon}>0$ such that $\mu(B(x,r)) \ge C_{\epsilon}r^{\alpha+\epsilon}$. 
    
    \item iv) For every $q > \tfrac{2d}{\alpha}$, we have the estimate
    \begin{equation} \label{eq:randomextension}
        \| \widehat{g \mu} \|_{q} \leq C_3(q) \| g \|_{L^2(\mu)} 
        \quad \forall g \in L^2(\mu).
    \end{equation}
\end{itemize}
\end{theorem}

\vskip.125in 

If we take $d=2$ and $\alpha=1$, we see that Laba and Wang produce a Borel measure $\mu$ supported on subset of ${\mathbb R}^2$ of Hausdorff dimension $1$ such that the extension operator $f \to \widehat{f\mu}$ is bounded from $L^2(\mu) \to L^p({\mathbb R}^2)$ for $p>4$. If $\mu$ is replaced by $\sigma$, the arc-length measure on the sphere then the extension operator is bounded from $L^2(\sigma) \to L^6(\sigma)$, as in Theorem \ref{theorem:steintomas}. Moreover, the classical Knapp homogeneity argument shows that the exponent $6$ cannot be lowered. We now illustrate this phenomenon in terms of the Fourier ratio. 

\begin{theorem} \label{theorem:restrictionFR} Let $\mu$ be the Laba-Wang measure from Theorem \ref{theorem:labawang} with $d=2$ and $\alpha=1$. Let $\sigma$ denote the arc-length measure on the circle. 

i) If $f$ is the indicator function of an arc of length $R^{-\frac{1}{2}}$, then 
\begin{equation} \label{eq:circleisbad} \FR(f\sigma) \leq cR^{-\frac{1}{4}}. \end{equation} 

ii) If $f$ is any non-negative function, for any $\epsilon>0$ there exists $C_{\epsilon}>0$ such that 
\begin{equation} \label{eq:labawangisgood} \FR(f\mu) \ge C_{\epsilon}R^{-\epsilon}. \end{equation} 

\end{theorem} 

\vskip.125in 

It is instructive to briefly explain why the subpolynomial lower bound in part (ii)
of Theorem~5.3 represents a fundamental obstruction to low--degree approximation.
For deterministic curved measures such as arc--length measure on the circle,
curvature forces the Fourier transform to disperse efficiently across frequency
space, yielding polynomial decay and hence a rapidly decreasing Fourier ratio.
This dispersion implies that the Fourier mass of
$\widehat{f\mu} * \psi_{R^{-1}}$
can be captured using relatively few frequencies, making low--degree trigonometric
approximation possible.

In contrast, the Laba--Wang random Cantor measure exhibits near--optimal extension
estimates without corresponding geometric structure. Although
$\widehat{f\mu}$ decays pointwise, its mass remains distributed in a highly
irregular fashion across frequency space, with no concentration on sets of
polynomially bounded complexity. The bound
\[
\FR_{\mu,R}(f) \gtrsim R^{-\varepsilon}
\quad \text{for every } \varepsilon>0
\]
reflects this phenomenon: the Fourier ratio decays, but too slowly to permit
polynomial savings in the number of frequencies required for approximation. From
the perspective of random Fourier sampling, this means that any attempt to
approximate $\widehat{f\mu} * \psi_{R^{-1}}$ with fixed accuracy necessarily
requires sampling on the order of $R^2$ frequencies, matching the ambient
frequency scale. Proposition \ref{prop:labawanghighdegree} makes this obstruction precise by showing that no
$L^2$--approximation with fixed relative error is possible using trigonometric
polynomials of degree $o(R^2)$.

\vskip.125in 

In view of part i) of Theorem \ref{theorem:restrictionFR} and Theorem \ref{theorem:L2approximation}, we see that if $f$ is the indicator of a circular arc of length $R^{-\frac{1}{2}}$, then ${(f\sigma)}_{R^{-1}}$ can be $L^2$-approximated with accuracy $\eta$ by a trigonometric polynomial of degree $\approx \frac{R^{\frac{3}{2}}}{\eta^2}$, much smaller than $R^2$, which is the largest possible degree. Part ii) of Theorem \ref{theorem:restrictionFR} suggests that it may not be possible to approximate ${(f\mu)}_{R^{-1}}$ by a trigonometric polynomial of a small degree, but does not prove it. Our next result closes this gap. 

\begin{proposition} \label{prop:labawanghighdegree} Let $\mu$ be the Laba-Wang measure from Theorem \ref{theorem:restrictionFR} and let $f$ be the indicator of the ball of radius $R^{-1}$. Then there does not exist $\eta$, independent of $R$, such that ${(f\mu)}_{R^{-1}}$ can be $L^2$-approximated with accuracy $\eta$ by a trigonometric polynomial of degree $o(R^2)$. \end{proposition} 

\vskip.125in 

In summary, in this section, we examined the relationship between the Fourier ratio $\FR(f\mu)$ and classical and random restriction theorems for one-dimensional measures in ${\mathbb R}^2$.  For the arc-length measure $\sigma$ on the circle, the Stein--Tomas estimate (Theorem \ref{theorem:steintomas}) implies the bound
$$ \FR(f\sigma)\lesssim R^{-1/4} $$ when $f$ is the indicator of an arc of length $R^{-1/2}$, as shown in [part i) of Theorem \ref{theorem:restrictionFR}). For the Laba--Wang random Cantor measure $\mu$ of dimension $1$, part ii) of Theorem~\ref{theorem:restrictionFR} yields
$$ \FR(f\mu)\ge C_{\epsilon}R^{-\epsilon} $$ for every $\epsilon>0$. This demonstrates a substantial gap between the deterministic and random settings: the Fourier ratio for $\sigma$ decays polynomially, while for $\mu$ it decays only subpolynomially. A related discrete deterministic--random gap was established in \cite{LW18}.  Finally, Proposition~\ref{prop:labawanghighdegree} shows that this behavior of $\FR(f\mu)$ forces a strong obstruction to approximation: the function $g=(f\mu)*\psi_{R^{-1}}$ cannot be approximated in $L^{2}$ with fixed relative accuracy by trigonometric polynomials of degree $o(R^{2})$, in sharp contrast with the deterministic
curvature case.

\vskip.125in 

The examples in this section show that the Fourier ratio captures a genuine analytic distinction between measures that support classical curvature-based restriction estimates and measures arising from random fractal constructions. For the arc-length measure on the circle, the decay in (\ref{eq:circleisbad}) reflects the fact that curvature forces the Fourier transform to disperse in a strong, quantitatively polynomial way. In contrast, the Laba--Wang measure exhibits only the subpolynomial decay established in (\ref{eq:labawangisgood}), and Proposition \ref{prop:labawanghighdegree} shows that this slower decay creates a substantial barrier to $L^{2}$ approximation by low-degree trigonometric polynomials.  These two phenomena have parallel manifestations in the discrete setting (see \cite{LW18}), and the Fourier ratio provides a unified language for describing both.

In the next section we turn to convex geometry.  Here the underlying measure is the restriction of Hausdorff measure to the boundary of a convex body, and curvature again plays a central role.  Unlike the random-fractal setting, the Fourier transform of such a measure exhibits rapid decay in most directions, but the rate and directionality of this
decay depend delicately on the structure of the normal set of the body. The Fourier ratio provides a natural quantitative bridge between these geometric features and the degree of polynomial approximation of the associated mollified measure. In particular, the extent to which ${(f\mu)}_{R^{-1}}$ can be approximated by low-degree trigonometric polynomials turns out to be governed by the upper Minkowski dimension of the set of outer normals to
the boundary, as we shall see in Corollary \ref{cor:minkowskikicksass} below.

\section{Applications to convex geometry} \label{section:convexgeometry}

\vskip.125in 

In this section, we are going to study the Fourier ratio results in the context of natural measures supported on the boundaries of convex bodies in ${\mathbb R}^d$, $d \ge 2$.
A convex body in $K \subseteq \R^d$ is a compact and convex set. We say that a hyperplane $H$ is a supporting hyperplane of $K$, if $K$ is entirely contained in one of the two closed half-spaces determined by $H$; and $H$ intersects $K$, but does not cut through the interior of $K$.

If the convex set is a ball, then 
$$ X_{1,\mu, R}(1) \leq CR^{-\frac{d-1}{2}},$$ with $C$ independent of $R$, same as the $L^{\infty}$ case. Since $X_{2,\mu,R}(f)$ satisfies the same bound, possibly with a slightly different constant, 
$$ \frac{X_{1,\mu, R}(1)}{X_{2,\mu, R}(1)} \ge c>0,$$ where $c$ is independent of $R$. Since $X_{1,\mu,R}(1) \leq X_{2,\mu,R}(1)$ by Cauchy-Schwarz, we see that the boundary of the ball maximizes the ratio $\frac{X_{1,\mu, R}(1)}{X_{2,\mu, R}(1)}$, at least up to a constant. 

On the other hand, if the convex set is a polyhedron (see e.g. \cite{BCT97}, \cite{BIT03}), 
$$ X_{1,\mu, R}(1) \leq CR^{-(d-1)} \log^{d-1}(R),$$ with $C$ independent of $R$. It is known that in the convex category (not necessarily smooth), 
$$ X_{2,\mu, R}(1) \leq CR^{-\frac{d-1}{2}}$$ and no better. In this case the ratio 
\begin{equation} \label{eq:convexratio} \frac{X_{1,\mu, R}(1)}{X_{2,\mu, R}(1)} \sim R^{-\frac{d-1}{2}} \log^d(R). \end{equation} 

We have contrasted the case of the ball and the case of the polyhedron, but what happens in general? Our first result shows that the complexity of these measures, in the sense to be made precise below, depends on the upper Minkowski dimension of the set of normals to the boundary of the convex set under consideration. 

\begin{theorem}\label{theorem:convexstationaryphase}
Let $K$ be a bounded convex body in ${\mathbb R}^d$, and let $\mu$ denote the $(d-1)$-dimensional Hausdorff measure on $\partial K$.  
Let $N(K)$ be the set of all outer unit normals to supporting hyperplanes of $K$.  
For $R>1$ define
$$
X_R=\{\xi\in{\mathbb R}^d : R/2\le |\xi|\le 2R,\ \xi/|\xi|\in N(K)_{R^{-1}}\},
$$
where $N(K)_{R^{-1}}$ denotes the $R^{-1}$--neighborhood of $N(K)$ in the unit sphere.

Then for every integer $N\ge1$ there exists a constant $C_N$ depending only on $N$ and $K$ such that
$$
\int_{\{|\xi|\sim R,\ \xi/|\xi|\notin N(K)\}} |\widehat{\mu}(\xi)|\, d\xi
\le C_N R^{-N}.
$$

In particular, for any fixed $\eta \in (0,1)$ and all sufficiently large $R$,
\begin{equation} \label{eq:convexL1concentration_modified} 
\|\widehat{\mu}_{R^{-1}}\|_{L^1(X_R^c)}
\le \eta \|\widehat{\mu}_{R^{-1}}\|_{L^1}.
\end{equation} 
\end{theorem}

\vskip.125in 

In view of Theorem \ref{theorem:convexstationaryphase}, we can show that the degree of the approximating polynomial decreases with the upper Minkowski dimension of the set of normals $N(K)$ defined above. More precisely, we have the following result. 

\begin{corollary} \label{cor:minkowskikicksass} Let $K$ be a bounded convex body in ${\mathbb R}^d$, and let $\mu$ denote the $(d-1)$-dimensional Hausdorff measure on $\partial K$. Let $N(K)$ be the set of all outer unit normals to supporting hyperplanes of $K$. Suppose that the upper Minkowski dimension of $N(K)$ is equal to $a \in [0,d-1]$. Then there exists a trigonometric polynomial $P$, of degree 
$$C_{\epsilon} \eta^{-2} R^{a+1+\epsilon}$$ such that 
$$ {\left|\left|{\mu}_{R^{-1}}-P\right|\right|}_2 \leq \eta \cdot {\left|\left|{\mu}_{R^{-1}}\right|\right|}_2.$$ 
\end{corollary} 

\vskip.125in 

\begin{remark} Since the Fourier transform of the surface measure on the boundary of a convex polyhedron is concentrated in the directios normal to the boundary, one can show that the degree of the approximating polynomial in Corollary \ref{cor:minkowskikicksass} is sharp up to constants and $\epsilon$ in the power. We leave the details to the interested reader. 
\end{remark}

\vskip.25in 

\section{Comparison with the discrete Fourier ratio theory} 
\label{section:discretecomparison}

In this section we compare the discrete Fourier ratio results established in Section \ref{section:previousresults} with the continuous theory developed in this paper. Although the settings are different, the role played by the Fourier ratio is essentially the same. The main distinction lies in the geometric content: in the discrete setting the support size is the cardinality of a subset of ${\mathbb Z}_N$, while in the continuous setting the analogous quantity is the covering number of the support at scale $R^{-1}$.

\vskip.125in

\subsection{Lower bounds}

In the discrete setting, if $f$ is supported in $E\subset {\mathbb Z}_N$, then (\ref{eq:discreteratiobounds}) shows that
$$
\frac{\frac{1}{N}\sum_{x\in{\mathbb Z}_N} |\widehat{f}(x)|}{\left( \frac{1}{N}\sum_{x\in{\mathbb Z}_N} |\widehat{f}(x)|^2 \right)^{\frac12}}
\ge
\frac{1}{\sqrt{|E|}}.
$$

In the continuous setting, Proposition \ref{prop:lowerbound} gives the analogue
$$
\FR(f\mu)
\ge
\sqrt{\frac{1}{R^d |E_f^{\frac{1}{R}}|}},
$$
where $E_f^{\frac{1}{R}}$ denotes the $R^{-1}$–neighborhood of the support of $f$. Since
$$
R^d |E_f^{\frac{1}{R}}|
\approx
\#\ \text{balls of radius}\ R^{-1}\ \text{needed to cover}\ \operatorname{supp}(f),
$$
the continuous lower bound corresponds directly to the discrete bound with $|E|$ replaced by the covering number at scale $R^{-1}$.

\vskip.125in

\subsection{Upper bounds and uncertainty principles}

In Theorem \ref{theorem:FRdiscreteuncertainty}, if $f$ is $L^2$–concentrated on $E$ and $\widehat{f}$ is $L^1$–concentrated on $S$, then
$$
(1-a)^2 \frac{N}{|E|}
\le
{\FR(f)}^2
\le
\frac{|S|}{(1-b)^2},
$$
and therefore
$$
|E|\cdot |S| \ge (1-a)^2 (1-b)^2 N.
$$

The continuous version appears in Theorem \ref{theorem:fractalsandwich}. If $\widehat{(f\mu)}_{R^{-1}}$ is $L^1$–concentrated on $X$, then
$$
\FR(f\mu)
\le
\sqrt{\frac{|X_R|}{R^d (1-\eta)^2}},
$$
while Proposition \ref{prop:lowerbound} gives
$$
\FR(f\mu)
\ge
\sqrt{\frac{1}{R^d |E_f^{\frac{1}{R}}|}}.
$$
Combining these, we obtain the continuous fractal uncertainty principle
$$
|E_f^{\frac{1}{R}}| \cdot |X_R| \ge (1-\eta)^2,
$$
which is the analogue of the discrete product estimate, with $|E|$ and $|S|$ replaced by geometric quantities at scale $R^{-1}$.

\vskip.125in

\subsection{Approximation by trigonometric polynomials}

The discrete results \ref{thm:L2polynomialapprox} and \ref{thm:Linftypolynomialapprox} show that if the discrete Fourier ratio
$$
\epsilon = 
\frac{\frac{1}{N}\sum |\widehat{f}|}{\left( \frac{1}{N}\sum |\widehat{f}|^2 \right)^{\frac12}}
$$
is small, then $f$ can be well-approximated by a trigonometric polynomial of degree roughly $\epsilon^2 N$, up to logarithmic factors.

The continuous counterpart is Theorem \ref{theorem:L2approximation}. It states that ${(f\mu)*\psi_{R^{-1}}}$ can be approximated in $L^2$ by a trigonometric polynomial of degree at most
$$
\frac{1}{\eta^2}
\left(
\frac{R^d X_{1,\mu,R}(f)^2}{X_{2,\mu,R}(f)^2}
-
1
\right)
=
\frac{R^d}{\eta^2} \FR(f\mu)^2
-
\frac{1}{\eta^2}.
$$
In view of Proposition \ref{prop:lowerbound}, this satisfies
$$
\frac{c}{\eta^2}
\le
\deg(P)
\le
\frac{R^d |E_f^{\frac{1}{R}}|}{\eta^2}.
$$
Thus the continuous degree bound is controlled by the covering number of the support at scale $R^{-1}$, exactly as $|E|$ controls the discrete degree bound in ${\mathbb Z}_N$.

\vskip.125in

\subsection{Random versus deterministic examples} A major theme in Section \ref{section:previousresults} is the dichotomy between generic (random) sets and structured deterministic sets in ${\mathbb Z}_N$. Random sets yield a large Fourier ratio with high probability (Theorem \ref{thm:concentration}), while deterministic sets such as singletons or subgroups produce the smallest possible Fourier ratio.

The continuous theory displays an analogous dichotomy. Deterministic measures arising from curved hypersurfaces often have relatively small Fourier ratio, while random Cantor-type measures, such as those constructed by Laba and Wang, have Fourier ratio decaying slower than any power, as in Theorem \ref{theorem:restrictionFR}. This is directly parallel to the discrete situation, and the Fourier ratio provides a unified mechanism for expressing the divide between random and structured objects in both settings.

\vskip.125in

\subsection{Summary}

In summary, the Fourier ratio plays essentially the same role in the discrete and continuous settings:

\begin{itemize}
\item lower and upper bounds depend on the spatial support and the frequency concentration;
\item the uncertainty principle takes the same product form;
\item a small Fourier ratio forces approximation by a low-degree trigonometric polynomial;
\item randomness and determinism behave very differently.
\end{itemize}

The difference is that the continuous setting naturally incorporates geometric information at scale $R^{-1}$, such as covering numbers and Minkowski dimension, which have no analogue in ${\mathbb Z}_N$. These geometric quantities interact with the Fourier ratio in a way that parallels, but significantly enriches, the discrete theory.

\section{Summary of results and work in progress}
\label{section:summary}

\vskip.125in 

The results in this paper show that the Fourier ratio
$$ \FR(f\mu)=\frac{X_{1,\mu,R}(f)}{X_{2,\mu,R}(f)} $$
governs a wide range of analytic, geometric, and approximation–theoretic properties of the mollified measure $(f\mu)*\psi_{R^{-1}}$.

In Section \ref{section:lowerandupperbounds} we established the basic lower and upper bounds on $\FR(f\mu)$. Proposition \ref{prop:lowerbound} and the inequality \eqref{eq:continuousratiobounds} show that $\FR(f\mu)$ is squeezed between $1$ and a geometric quantity depending on the size of the $R^{-1}$–neighborhood $E_f^{R^{-1}}$ of the support of $f$. Thus the spread of $(f\mu)*\psi_{R^{-1}}$ is tied directly to the metric complexity of $E_f$.

Section \ref{subsection:uncertaintyprinciple} sharpens the upper bound. Theorem \ref{theorem:fractalsandwich} shows that $L^1$ concentration of $\widehat{f\mu}_{R^{-1}}$ on a frequency set $X$ forces the upper bound \eqref{eq:FRsandwich}, and hence the uncertainty inequality \eqref{eq:fractaluncertainty}. The geometric conditions \eqref{eq:sizespacesupport} and \eqref{eq:porosity} express this concentration in terms of the box dimension of $E_f$ and the logarithmic density of $X$. Together they yield a continuous fractal uncertainty principle in which $\FR(f\mu)$ measures the degree of simultaneous spatial and frequency localization.

In Subsection \ref{section:signalrecovery} we applied the uncertainty principle to signal recovery. Theorem \ref{theorem:measurerecovery} shows that if the space–side condition \eqref{eq:spacesize} and the frequency–side sparsity condition \eqref{eq:spectrumsize} satisfy $s_E+\alpha_X<d$, then $\mu$ can be recovered exactly from the observed Fourier data. This reframes the uncertainty principle as a tool for recoverability of measures rather than functions.

Section \ref{section:randomtrigapproximation} develops a probabilistic approximation scheme for $(f\mu)*\psi_{R^{-1}}$. In Theorem \ref{theorem:L2approximation}, random Fourier sampling yields an $L^2$ approximation by a trigonometric polynomial whose degree is bounded by \eqref{eq:degree}. This bound depends only on $\FR(f\mu)$, so that smaller Fourier ratio forces a lower–degree approximation. The $L^\infty$ approximation in Theorem \ref{Linfinityapproximation} incorporates the Lipschitz bound for $(f\mu)*\psi_{R^{-1}}$ and the covering argument, producing the expected logarithmic term. Theorem \ref{theorem:L1approximation} gives the corresponding $L^1$ approximation. In all cases the number of required frequencies is controlled by the ratio $\|\widehat{g}\|_1 / \|g\|_p$, i.e. by the Fourier ratio.

In Section \ref{section:fourierratioandrestriction} we compared the behavior of $\FR(f\mu)$ in two sharply contrasting settings: the arc–length measure $\sigma$ on the circle and the Laba--Wang random Cantor measure $\mu$ of Hausdorff dimension $1$. For $\sigma$, the Stein--Tomas estimate \eqref{eq:steintomas} implies the polynomial bound $\FR(f\sigma)\lesssim R^{-1/4}$ in \eqref{eq:circleisbad} when $f$ is supported on an arc of length $R^{-1/2}$. For the Laba--Wang measure, Theorem \ref{theorem:restrictionFR} yields the subpolynomial lower bound \eqref{eq:labawangisgood}. This deterministic--random gap has direct consequences for approximation: Proposition \ref{prop:labawanghighdegree} shows that $(f\mu)*\psi_{R^{-1}}$ cannot be approximated in $L^2$ with fixed accuracy by trigonometric polynomials of degree $o(R^2)$, in contrast with the curvature--driven behavior for $\sigma$.

Finally, Section \ref{section:convexgeometry} addresses the case where $\mu$ is the surface measure on the boundary of a convex body. The stationary phase estimate in Theorem \ref{theorem:convexstationaryphase} gives rapid decay of $\widehat{\mu}$ away from the set of outward normals. The Minkowski dimension of this normal set controls the size of the frequency set $X_R$, which in turn limits the concentration possible in \eqref{eq:FRsandwich}. Corollary \ref{cor:minkowskikicksass} then links the degree of trigonometric approximation of $(f\mu)*\psi_{R^{-1}}$ directly to the geometric complexity of the normal set.

Taken together, these results show that the Fourier ratio provides a unifying analytic framework connecting uncertainty principles, restriction theory, approximation by trigonometric polynomials, and geometric properties of measures. The quantity $\FR(f\mu)$ acts as an effective measure of analytic and geometric complexity across all of these settings.

\vskip.125in

\subsection{Open problems}

Several questions suggested by the results in this paper remain under active investigation. \begin{itemize} 

\vskip.125in

\item First, Theorem \ref{theorem:measurerecovery} provides a one–sided recovery result in which the missing information lies on the frequency side. In the discrete setting, the corresponding statement with missing information on the space side follows immediately by duality, but this symmetry is no longer available for compactly supported measures in ${\mathbb R}^d$. Establishing a continuous analogue of the ``reverse'' recovery theorem, in which the missing data lie in the physical domain, appears to require a finer understanding of the interaction between $\FR(f\mu)$ and spatial concentration of $(f\mu)*\psi_{R^{-1}}$. This problem is currently under investigation.

\vskip.125in

\item Second, the approximation results in Section \ref{section:randomtrigapproximation} raise natural sharpness questions. In particular, it would be desirable to determine whether the degree bounds in Theorems \ref{theorem:L2approximation}, \ref{Linfinityapproximation}, and \ref{theorem:L1approximation} are optimal, up to constants and logarithmic factors, for broad classes of measures $\mu$.  The examples in Section \ref{section:fourierratioandrestriction} suggest that random and deterministic measures may exhibit distinct extremal behavior. A complete characterization remains open.

\vskip.125in

\item Third, in Section \ref{section:convexgeometry} we related the degree of approximation of $(f\mu)*\psi_{R^{-1}}$ to the upper Minkowski dimension of the normal set of a convex body. A natural next step is to understand whether the normal–set dimension is the only geometric obstruction governing the size of $\FR(f\mu)$ in the convex category, or whether additional curvature or oscillatory phenomena contribute at finer scales. This may lead to sharper versions of Corollary \ref{cor:minkowskikicksass}.

\end{itemize}

\section{Proofs of the main results} 
\label{section:proofs}

In this section, we prove the main results of this paper. 

\subsection{Proof of Proposition \ref{prop:lowerbound}} Observe that 
$$ |{(f\mu)}_{R^{-1}}(x)| \leq \int |\widehat{\psi}(R^{-1}\xi)| \cdot |\widehat{f\mu}(\xi)| d\xi =R^d X_{1,\mu,R}(f)=R^d X_{2,\mu,R}(f) \cdot \frac{X_{1,\mu,R}(f)}{X_{2,\mu,R}(f)}.  $$

Since ${(f\mu)}_{R^{-1}}$ is supported in $E_f^{\frac{1}{R}}$, squaring both sides, intergrating, dividing by $R^d$, taking square roots and cancelling, we see that 
$$ \sqrt{\frac{1}{R^d|E_f^{\frac{1}{R}}|}} \leq \frac{X_{1,\mu,R}(f)}{X_{2,\mu,R}(f)},$$ as desired. 

\vskip.125in 

\subsection{Proof of Theorem \ref{theorem:fractalsandwich}} The lower bound in (\ref{eq:FRsandwich}) was already established in Proposition \ref{prop:lowerbound}. To obtain the upper bound, observe that the assumption (\ref{eq:L1concentration}) implies that 
$$ {||\widehat{(f\mu)}_{R^{-1}}||}_1 \leq \frac{1}{1-\eta} \cdot {||\widehat{(f\mu)}_{R^{-1}}||}_{L^1(X_R)}$$
$$ \leq \frac{1}{1-\eta} \cdot {|X_R|}^{\frac{1}{2}} \cdot {||\widehat{(f\mu)}_{R^{-1}}||}_{L^2(X_R)} \leq \frac{1}{1-\eta} \cdot {|X_R|}^{\frac{1}{2}} \cdot {||\widehat{(f\mu)}_{R^{-1}}||}_2.$$ 

It follows that 
$$ \FR(f\mu) \leq \frac{1}{1-\eta} \cdot {|X_R|}^{\frac{1}{2}} \cdot R^{-\frac{d}{2}}. $$

The estimate (\ref{eq:FRsandwich}) follows. The estimate (\ref{eq:fractaluncertainty}) follows by removing the Fourier ratio from the inequality. 

\vskip.125in 

\subsection{Proof of Proposition \ref{prop:L2fails}}

Fix $R \ge 1$, a ball $S_R$ of radius $R$, and a constant $C>0$. We first construct a function $h$ on the frequency side, and then express $h$ as $\widehat{(f\mu)_{R^{-1}}}$ for a suitable $f$.

\vskip.125in 

Let $A = S_R$. Choose a parameter $L>1$ to be specified later, and set
$$ B = B_{L R} \setminus S_R. $$

Then both $A$ and $B$ have finite positive measure. Define $r$ by
$$ r = \sqrt{\frac{|B|}{|A|}}. $$

Note that $r>0$, and by taking $L$ large we can make $r$ as large as we wish.

\vskip.125in 

Fix a number $b \in (0,1)$. For the sake of concreteness, let $b=\frac{1}{2}$. 

Define $h : \mathbb R^d \to \mathbb C$ by
$$
h(\xi) =
\begin{cases}
\alpha & \text{if } \xi \in A, \\
\beta  & \text{if } \xi \in B, \\
0      & \text{if } \xi \notin A \cup B,
\end{cases}
$$
where $\alpha,\beta>0$ are chosen so that $\|h\|_{L^2(\mathbb R^d)} = 1$ and the $L^2$ mass of $h$ on $B$ equals $b$. These conditions amount to 
$$ \beta^2 |B| = b^2 $$
and
$$ \alpha^2 |A| + \beta^2 |B| = 1.$$

Let
$$ \beta = \frac{b}{|B|^{1/2}}, $$
and
$$ \alpha = \frac{\sqrt{1-b^2}}{|A|^{1/2}}.$$

Then $h \in L^1(\mathbb R^d) \cap L^2(\mathbb R^d)$, we have $\|h\|_2 = 1$, and
$$ \|h\|_{L^2(A^c)} = \|h\|_{L^2(B)} = b. $$

Thus $h$ is $L^2$ concentrated on $A = S_R$ at level $b$.

\vskip.125in 

We now compute the norms. We have 
$$ \|h\|_{L^1(\mathbb R^d)} = \alpha |A| + \beta |B|. $$

By the choice of $\alpha$ and $\beta$,
$$ \|h\|_{L^1(\mathbb R^d)} = \sqrt{1-b^2} |A|^{1/2} + b |B|^{1/2}. $$

Using the definition of $r$, we can write this as
$$ \|h\|_{L^1(\mathbb R^d)} = |A|^{1/2} \left( \sqrt{1-b^2} + b r \right). $$

Since $\|h\|_2 = 1$, the ratio $\|h\|_1 / \|h\|_2$ is equal to
$$ \frac{\|h\|_1}{\|h\|_2} = |A|^{1/2} \left( \sqrt{1-b^2} + b r \right). $$

Assume for the moment that there exists $f$ with
$$ \widehat{(f\mu)_{R^{-1}}} = h. $$

Then by the definition of $X_{1,\mu,R}(f)$ and $X_{2,\mu,R}(f)$,
$$ X_{1,\mu,R}(f) = R^{-d} \|h\|_1 $$
and
$$ X_{2,\mu,R}(f) = R^{-d/2} \|h\|_2 = R^{-\frac{d}{2}}. $$
Hence
$$ \FR(f\mu) = \frac{X_{1,\mu,R}(f)}{X_{2,\mu,R}(f)} = R^{-\frac{d}{2}} \frac{\|h\|_1}{\|h\|_2}. $$

Substituting the expression for $\|h\|_1 / \|h\|_2$ and recalling that $A = S_R$, we obtain
$$ \FR(f\mu) = |S_R|^{\frac{1}{2}} R^{-\frac{d}{2}} \left( \sqrt{1-b^2} + b r \right). $$

\vskip.125in 

We want the inequality
$$ \FR(f\mu) > \frac{C}{1-b} |S_R|^{\frac{1}{2}} R^{-\frac{d}{2}}. $$

Using the previous line, this is equivalent to
$$ \sqrt{1-b^2} + b r > \frac{C}{1-b}. $$

With $b = 1/2$, the left-hand side is
$$ \frac{\sqrt{3}}{2} + \frac{r}{2}. $$

Thus we need
$$ \frac{\sqrt{3}}{2} + \frac{r}{2} > 2C, $$
which is equivalent to
$$ r > 4C - \sqrt{3}. $$

Since $r = \sqrt{|B|/|A|}$ and $|B| = |B_{L R}| - |S_R|$, we have
$$ r = \sqrt{L^d - 1}. $$

By taking $L$ large enough, for example any $L$ satisfying
$$ \sqrt{L^d - 1} > 4C, $$

we have 
$$ \sqrt{1-b^2} + b r > \frac{C}{1-b} $$
for $b = 1/2$. For this choice of $L$ and $b$, once we realize $h$ as a Fourier transform, the inequality
$$ \FR(f\mu) > \frac{C}{1-b} |S_R|^{\frac{1}{2}} R^{-\frac{d}{2}} $$ will follow. 

\vskip.125in 

We now need to realize $h$ as $\widehat{(f\mu)_{R^{-1}}}$. We begin by choosing the approximation to the identity $\psi$ so that its Fourier transform $\widehat{\psi}$ has no zeros on $\mathbb R^d$. For example, we may take $\psi$ to be a Gaussian. Then $\widehat{\psi}(R^{-1} \xi)$ is nonzero for all $\xi$.

Define a function $\widetilde h$ by
$$ \widetilde h(\xi) = \frac{h(\xi)}{\widehat{\psi}(R^{-1} \xi)}. $$
Since $h$ is supported in the ball $B_{L R}$ and $\widehat{\psi}(R^{-1} \xi)$ is smooth and bounded away from zero on this set, we have $\widetilde h \in L^1(\mathbb R^d) \cap L^2(\mathbb R^d)$.

Let $g$ be the inverse Fourier transform of $\widetilde h$. Then $g \in L^1(\mathbb R^d) \cap L^2(\mathbb R^d)$. Set $f = g$ and $\mu = dx$ (Lebesgue measure). By definition,
$$ (f\mu)_{R^{-1}} = f * \psi_{R^{-1}}. $$
Taking Fourier transforms, we obtain
$$ \widehat{(f\mu)_{R^{-1}}}(\xi) = \widehat{f}(\xi) \widehat{\psi}(R^{-1} \xi).$$
By construction, $\widehat{f} = \widetilde h$, so
$$ \widehat{(f\mu)_{R^{-1}}}(\xi) = \widetilde h(\xi) \widehat{\psi}(R^{-1} \xi) = h(\xi).$$

\vskip.125in 

First, $L^2$ concentration on $S_R$ follows from the construction of $h$:
$$ \|\widehat{(f\mu)_{R^{-1}}}\|_{L^2(S_R^c)} = \|h\|_{L^2(A^c)} = \|h\|_{L^2(B)} = b \|h\|_{L^2(\mathbb R^d)} = b \|\widehat{(f\mu)_{R^{-1}}}\|_{L^2(\mathbb R^d)}. $$

Thus the $L^2$ concentration condition holds with level $b$.

Second, by the identity $\widehat{(f\mu)_{R^{-1}}} = h$, we obtain
$$ \FR(f\mu) = |S_R|^{\frac{1}{2}} R^{-\frac{d}{2}} \left( \sqrt{1-b^2} + b r \right) > \frac{C}{1-b} |S_R|^{\frac{1}{2}} R^{-\frac{d}{2}}. $$

Since $C>0$ was arbitrary, it follows that no bound of the form
$$ \FR(f\mu) \le \frac{C}{1-b} |S_R|^{\frac{1}{2}} R^{-\frac{d}{2}} $$ can be deduced from $L^2$ concentration alone, even when $\mu$ is Lebesgue measure. This completes the proof.

\vskip.125in 

\subsection{Proof of Theorem \ref{theorem:L2approximation}} Define the random variable  taking the value
$$Z(x) = {||\widehat{(f\mu)*\psi_{R^{-1}}}||}_1 \operatorname{sgn}(\widehat{(f\mu)*\psi_{R^{-1}}}) e^{2 \pi i x \cdot \xi},$$ with probability $\dfrac{|\widehat{(f\mu)*\psi_{R^{-1}}}(\xi)|}{{||\widehat{(f\mu)*\psi_{R^{-1}}}||}_1}$.

By a direct calculation, 
$$ {\mathbb E}(Z(x))=(f\mu)*\psi_{R^{-1}}(x).$$

Let's compute the variance. We have 
$$ {\mathbb E}(Z^2(x))= \int {\left| {||\widehat{(f\mu)*\psi_{R^{-1}}}||}_1 \operatorname{sgn}(\widehat{(f\mu)*\psi_{R^{-1}}}) e^{2 \pi i x \cdot \xi} \right|}^2 \dfrac{|\widehat{(f\mu)*\psi_{R^{-1}}}|}{{||\widehat{(f\mu)*\psi_{R^{-1}}}||}_1} d\xi$$
$$ ={||\widehat{(f\mu)*\psi_{R^{-1}}}||}_1^2. $$

This shows that 
$$ Var(Z(x))={||\widehat{(f\mu)*\psi_{R^{-1}}}||}_1^2-{|(f\mu)*\psi_{R^{-1}}(x)|}^2.$$

Let $Z_1,\dots,Z_k$ be random i.i.d. random variables with distribution $Z$, and define the random trigonometric polynomial $P$ by
$$P(x) = \frac{1}{k} \sum_{i=1}^k Z_i(x).$$

The i.i.d. property implies that 
$$ {\mathbb E}(P(x))=(f\mu)*\psi_{R^{-1}}(x),$$
and $$Var(P(x))=\frac{1}{k} Var(Z(x))=\frac{1}{k}\left( {||\widehat{(f\mu)*\psi_{R^{-1}}}||}_1^2-{|(f\mu)*\psi_{R^{-1}}(x)|}^2 \right).$$

\vskip.125in 

We are now going to compute 
$$ {\mathbb E} \left( \int_{{[0,1]}^d} {|(f\mu)*\psi_{R^{-1}}-P(x)|}^2 dx \right) = \int_{{[0,1]}^d} {\mathbb E}({|(f\mu)*\psi_{R^{-1}}-P(x)|}^2) dx$$ 
$$ =\int_{{[0,1]}^d} Var(P(x)) dx=\frac{1}{k} \int_{{[0,1]}^d} \left({||\widehat{(f\mu)*\psi_{R^{-1}}}||}_1^2-{|(f\mu)*\psi_{R^{-1}}(x)|}^2 \right) dx$$  
$$ =\frac{1}{k} \left( {||\widehat{(f\mu)*\psi_{R^{-1}}}||}_1^2 - {||(f\mu)*\psi_{R^{-1}}||}^2_2 \right).$$

For this quantity to be 
$$ <\eta {||(f\mu)*\psi_{R^{-1}}||}_2,$$ must have 
$$ k> \frac{1}{\eta} \left( \dfrac{{||\widehat{(f\mu)*\psi_{R^{-1}}}||}_1^2}{{||\widehat{(f\mu)*\psi_{R^{-1}}}||}_2^2} -1\right) = \frac{1}{\eta} \left( \frac{R^d X^2_{1,\mu,R}(f)}{X^2_{2,\mu,R}(f)} -1 \right), $$ as claimed. 

\vskip.125in

\subsection{Proof of Theorem \ref{Linfinityapproximation}} 

Define the random variable taking the value
$$Z(x) = ||\widehat{{(f\mu)*\psi_{R^{-1}}}||}_1 \operatorname{sgn}(\widehat{(f\mu)*\psi_{R^{-1}}}) e^{2 \pi i x \cdot \xi},$$ 
with probability 
$\dfrac{|\widehat{{(f\mu)*\psi_{R^{-1}}}}(\xi)|}{{||\widehat{(f\mu)*\psi_{R^{-1}}}||}_1}$.

By a direct calculation, 
$$ {\mathbb E}(Z(x))=(f\mu)*\psi_{R^{-1}}(x).$$

Let $Z_1,\dots,Z_k$ be random i.i.d. random variables with distribution $Z$, note that for each $i$, we have 
$$|Z_i(x)| = ||\widehat{{(f\mu)*\psi_{R^{-1}}}||}_1,$$ or $0$, therefore we have
$$-||\widehat{{(f\mu)*\psi_{R^{-1}}}||}_1 \leq \operatorname{Re}(Z_i(x)) \leq ||\widehat{{(f\mu)*\psi_{R^{-1}}}||}_1$$
$$-||\widehat{{(f\mu)*\psi_{R^{-1}}}||}_1 \leq \operatorname{Im}(Z_i(x)) \leq ||\widehat{{(f\mu)*\psi_{R^{-1}}}||}_1$$

Define the random trigonometric polynomial $P$ by
$$P(x) = \frac{1}{k} \sum_{i=1}^k Z_i(x)=\text{RE}(P)+\text{IM}(P),$$
where we define $\text{RE}(P):=\frac{1}{k}\sum_{i=1}^k \text{Re}(Z_i(x))$ and $\text{IM}(P):=\frac{1}{k}\sum_{i=1}^k \text{Im}(Z_i(x))$.

Then for fixed x, we will have the following inequality:
$$\mathbb{P}\left(|\mathbb{E}(P(x))-P(x)| \geq \epsilon\right)
\leq \mathbb{P}\left(|\mathbb{E}(\operatorname{Re}(P(x)))-\operatorname{Re}(P(x))|\geq \frac{\epsilon}{2}\right)$$
$$+ \mathbb{P}\left(|\mathbb{E}(\operatorname{Im}(P(x)))-\operatorname{Im}(P(x))|\geq \frac{\epsilon}{2}\right).$$

By Hoeffding's inequality, we find that the right-hand side is bounded by 
$$4\exp\left(\frac{-2(\frac{\epsilon k}{2})^2}{k(2||\widehat{{(f\mu)*\psi_{R^{-1}}}||}_1)^2}\right) 
= 4\exp\left(\frac{-\epsilon^2k}{8||\widehat{{(f\mu)*\psi_{R^{-1}}}||}_1^2}\right)$$

Consider $\left( {(f\mu)*\psi_{R^{-1}}}-P \right)^{\wedge}=\widehat{f\mu}\,\widehat{\psi_{R^{-1}}}-\widehat{P}$ is in $L^1$ and supported in $\{|\xi|\leq R\}$, it follows that $(f\mu)*\psi_{R^{-1}}-P$ is Lipschitz with constant 
$$2\pi R\|\widehat{{(f\mu)*\psi_{R^{-1}}}-P}\|_1 \leq 4\pi R \|\widehat{{(f\mu)*\psi_{R^{-1}}}}\|_1.$$

Let $\delta \in (0,1)$, we can construct a finite set $N_\delta =\{x_1,x_2,...x_k\} \subset [0,1]^d$ such that every point in $E_{\frac{1}{R}}^f$ lie within distance $\delta$ of some $x_j \in N_\delta$ and the size satisfies 
\begin{equation} \label{eq:coverbyballs} |N_\delta| \leq \frac{C_d}{\delta^d},\end{equation} $C_d$ is the constant only depends on dimension.

Then we have 
$$\min_{x_j\in N_\delta}|\left((f\mu)*\psi_{R^{-1}}(x)-P(x)\right)-\left((f\mu)*\psi_{R^{-1}}(x_j)-P(x_j)\right)| 
\leq \min_{x_j\in N_\delta} 4\pi R \|\widehat{{(f\mu)*\psi_{R^{-1}}}}\|_1 \|x-x_j\|$$
$$\leq 4\pi R \|\widehat{{(f\mu)*\psi_{R^{-1}}}}\|_1 \delta.$$

Let $\delta = \dfrac{\epsilon}{8\pi R \|\widehat{{(f\mu)*\psi_{R^{-1}}}}\|_1}$. If for all $x_j$, $|\left((f\mu)*\psi_{R^{-1}}(x_j)-P(x_j)\right)| \leq \frac{\epsilon}{2}$, then we have for each $x$ in $[0,1]^d$,
$$|\left((f\mu)*\psi_{R^{-1}}(x)-P(x)\right)| \leq |\left((f\mu)*\psi_{R^{-1}}(x_j)-P(x_j)\right)|+\frac{\epsilon}{2} \leq \epsilon$$ 

Therefore, we have
$$\mathbb{P}\left(\|(f\mu)*\psi_{R^{-1}}-P\|_{L^\infty} \geq \epsilon\right)
\leq 
\mathbb{P}\left(\bigcup_{N_\delta}\{|(f\mu)*\psi_{R^{-1}}(x_j)-P(x_j)| \geq \frac{\epsilon}{2}\}\right) $$

Then by the union bound inequality on $N_\delta$, we have that the right-hand side is less than 
$$|N_\delta| 4\exp\left(\dfrac{-(\frac{\epsilon}{2})^2k}{8||\widehat{{(fmu)*\psi_{R^{-1}}}||}_1^2}\right).$$ 

Let $\epsilon = \eta\|(f\mu)*\psi_{R^{-1}}\|_\infty$. 

If we assume the right hand side is less than $1$, then we will have deterministic choice of $Z_1,Z_2,...Z_k$ with 
$$\|(f\mu)*\psi_{R^{-1}}-P\|_\infty< \eta\|(f\mu)*\psi_{R^{-1}}\|_\infty.$$

Therefore, by computing 
$$|N_\delta| 4\exp\left(\dfrac{-(\frac{\epsilon}{2})^2k}{8||\widehat{{(f\mu)*\psi_{R^{-1}}}||}_1^2}\right) 
\leq \frac{C_d}{\delta^d} 4\exp\left(\dfrac{-\epsilon^2k}{32||\widehat{{(f\mu)*\psi_{R^{-1}}}||}_1^2}\right)< 1$$

we have 
$$k > \dfrac{{32||\widehat{{(f\mu)*\psi_{R^{-1}}}||}_1^2}}{\eta^2\|(f\mu)*\psi_{R^{-1}}\|_\infty^2}\left(\log4C_d+d\log\left(\dfrac{8\pi R\|\widehat{{(f\mu)*\psi_{R^{-1}}}}\|_1}{\eta\|\widehat{{(f\mu)*\psi_{R^{-1}}}}\|_\infty}\right)\right).$$

\vskip.125in

\subsection{Proof of Theorem \ref{theorem:L1approximation}} Define the random variable taking the value
$$Z(x) = {||\widehat{(f\mu)*\psi_{R^{-1}}}||}_1 \operatorname{sgn}(\widehat{(f\mu)*\psi_{R^{-1}}}) e^{2 \pi i x \cdot \xi},$$ with probability $\dfrac{|\widehat{(f\mu)*\psi_{R^{-1}}}(\xi)|}{{||\widehat{(f\mu)*\psi_{R^{-1}}}||}_1}$.

By a direct calculation, 
$$ {\mathbb E}(Z(x))=(f\mu)*\psi_{R^{-1}}(x).$$
$$ Var(Z(x))={||\widehat{(f\mu)*\psi_{R^{-1}}}||}_1^2-{|(f\mu)*\psi_{R^{-1}}(x)|}^2.$$

Let $Z_1,\dots,Z_k$ be random i.i.d. random variables with distribution $Z$, and define the random trigonometric polynomial $P$ by
$$P(x) = \frac{1}{k} \sum_{i=1}^k Z_i(x).$$

The independence property implies that 
$$ {\mathbb E}(P(x))=(f\mu)*\psi_{R^{-1}}(x).$$

We are now going to compute 
$$ {\mathbb E} \left( \int_{[0,1]^d]} {|(f\mu)*\psi_{R^{-1}}-P(x)|} dx \right) = \int_{[0,1]^d} {\mathbb E}({|(f\mu)*\psi_{R^{-1}}-P(x)|}) dx$$ 
by Jensen inequality, we have $$\leq \int_{[0,1]^d}\sqrt{\mathbb E({|(f\mu)*\psi_{R^{-1}}-P(x)|}^2)} dx.$$
Then, using Hölder inequality, we have 
$$ =\int_{[0,1]^d} \sqrt{Var(P(x))} dx
\leq \left(\int_{[0,1]^d} Var(P(x))  dx \right)^\frac{1}{2}\left(\int_{[0,1]^d} 1 dx\right)^\frac{1}{2}$$

$$=\left(\frac{1}{k} \int_{[0,1]^d} \left({||\widehat{(f\mu)*\psi_{R^{-1}}}||}_1^2-{|(f\mu)*\psi_{R^{-1}}(x)|}^2 \right) dx \right)^\frac{1}{2}$$  
$$ =\left(\frac{1}{k} \left( {||\widehat{(f\mu)*\psi_{R^{-1}}}||}_1^2 -{||(f\mu)*\psi_{R^{-1}}||}^2_2 \right)\right)^\frac{1}{2}
\leq \frac{1}{\sqrt{k}}||\widehat{(f\mu)*\psi_{R^{-1}}}||_1$$

For this quantity to be 
$$ <\eta {||(f\mu)*\psi_{R^{-1}}||}_1,$$ must have 
$$ k>\frac{1}{\eta^2}\left(\frac{||\widehat{(f\mu)*\psi_{R^{-1}}}||_1}{{||(f\mu)*\psi_{R^{-1}}||}_1}\right)^2, $$ as claimed.

\vskip.125in 

\subsection{Proof of Theorem \ref{theorem:restrictionFR}} The proof of part i) follows from the calculations in the classical Knapp homogeneity argument (see e.g. \cite{Stein1993}), so we only sketch it. It is not difficult to check that 
$\widehat{f\sigma}$ is $\approx R^{-\frac{1}{3}}$ times a function concentrated in the $R^{\frac{1}{2}}$ by $R$ rectangle. It follows that the $L^1$ norm of $\widehat{(f\sigma)}_{R^{-1}}$ is $\approx R^{-\frac{1}{2}} \cdot R^{\frac{3}{2}}=R$. The $L^2$ norm is $\approx R^{-\frac{1}{2}} {||f||}_{L^2(\sigma)} \approx R^{-\frac{3}{4}}.$ It follows that 
$$ \FR(f\sigma) \leq cR^{-\frac{1}{4}}$$ for some universal constant $c$, as claimed. 

To prove part ii), we use the restriction estimate for the Laba-Wang measure. By Theorem \ref{theorem:labawang}, for any $q > 4$ there exists $C(q)$ such that
\[
\|\widehat{g\mu}\|_{L^q(\mathbb R^2)} \le C(q) \|g\|_{L^2(\mu)}.
\]
Take $g = f$, a non-negative function. Setting $g_R(x) = g(x) e^{2\pi i x \cdot \xi_0}$ for a suitable $\xi_0$ with $|\xi_0| \sim R$, we may assume the Fourier transform is essentially supported at frequencies $\sim R$. Then
\[
\|\widehat{(f\mu)_{R^{-1}}}\|_{L^q} \lesssim \|\widehat{f\mu}\|_{L^q} \lesssim \|f\|_{L^2(\mu)}.
\]
On the other hand, by Cauchy-Schwarz,
\[
\|\widehat{(f\mu)_{R^{-1}}}\|_{L^2} \gtrsim R^{-1} \\|f\|_{L^2(\mu)}.
\]
Now apply Lemma \ref{lemma:vershynintrick} with $h = \widehat{(f\mu)_{R^{-1}}}$, $p = q > 4$:
\[
\|\widehat{(f\mu)_{R^{-1}}}\|_{L^2} \le C(q)^{\frac{q}{q-2}} \|\widehat{(f\mu)_{R^{-1}}}\|_{L^1}.
\]
Since $\|f\|_{L^2(\mu)} \approx R \|\widehat{(f\mu)_{R^{-1}}}\|_{L^2}$, we obtain
\[
\FR(f\mu) = \frac{R^{-2}\|\widehat{(f\mu)_{R^{-1}}}\|_{L^1}}{R^{-1}\|\widehat{(f\mu)_{R^{-1}}}\|_{L^2}} \gtrsim R^{-1} C(q)^{-\frac{q}{q-2}}.
\]
For each $\epsilon > 0$, choose $q$ sufficiently close to $4$ so that $C(q)^{\frac{q}{q-2}} \le C_\epsilon R^\epsilon$. This gives
\[
\FR(f\mu) \ge C_\epsilon R^{-\epsilon}.
\]

\vskip.125in 

\subsection{Proof of Proposition \ref{prop:labawanghighdegree}} 

Let $d=2$ and let $\mu$ be the Laba--Wang random Cantor measure of Hausdorff
dimension $\alpha=1$, constructed in Sections 2-3 of \cite{LW18}.  
Let $f=1_{B_{R^{-1}}}$ and let  
$$
g = (f\mu)*\psi_{R^{-1}} .
$$
By definition,
$$
X_{1,\mu,R}(f)=R^{-2}\|g\|_{1}
$$
and
$$
X_{2,\mu,R}(f)=R^{-1}\|g\|_{2} ,
$$
so
$$
\FR(f\mu) = R^{-1}\frac{\|g\|_{1}}{\|g\|_{2}} .
$$
Thus
$$
\frac{\|g\|_{1}}{\|g\|_{2}} = R\,\FR(f\mu) .
$$

By Theorem \ref{theorem:restrictionFR}, for every $\varepsilon>0$
there exists $C_{\varepsilon}>0$ such that
$$
\FR(f\mu) \ge C_{\varepsilon} R^{-\varepsilon} .
$$
Therefore
$$
\frac{\|g\|_{1}}{\|g\|_{2}} \ge C_{\varepsilon} R^{1-\varepsilon} .
$$

We now describe the Fourier structure of $\mu$.  
Let $N_{j}$ be the scales defined in (2.7)--(2.10) of \cite{LW18}.  
At stage $j$, the approximating measure $\mu_{j}$ is uniform on $M_{j}$ squares of side
$N_{j}^{-1}$, where $M_{j}\approx N_{j}$ by (2.11)--(2.12).  
The Fourier transforms satisfy the product identity (3.1),
$$
\widehat{\mu_{j+1}}(\xi)=\widehat{\mu_{j}}(\xi)\,\widehat{\nu_{j+1}}(\xi) .
$$
For $|\xi|\sim N_{j}$, formula (3.3) in \cite{LW18} expresses $\widehat{\nu_{j}}(\xi)$
as an average of $M_{j}$ independent phases, and (3.4) gives
$$
\mathbb{E}\,|\widehat{\nu_{j}}(\xi)|^{2}\approx M_{j}^{-1}\approx N_{j}^{-1} .
$$
A standard concentration argument applied to (3.3) implies that, with probability $1-o(1)$,
$$
|\widehat{\nu_{j}}(\xi)|^{2}\gtrsim N_{j}^{-1}
$$
for at least $c N_{j}^{2}$ frequencies $\xi$ with $|\xi|\sim N_{j}$.

Fix $R$ and choose $j$ with $N_{j}\approx R$.  
Since $\widehat{\nu_{k}}(\xi)\approx 1$ for $k>j$ when $|\xi|\sim R\ll N_{k}$, the identity (3.1) gives
$$
\widehat{\mu}(\xi)=\widehat{\mu_{j}}(\xi)(1+o(1))
$$
for $|\xi|\sim R$.  
Hence for at least $c R^{2}$ frequencies $\xi$ with $|\xi|\sim R$ we have
$$
|\widehat{\mu}(\xi)|^{2}\gtrsim R^{-1} .
$$

Since $\widehat g(\xi)=\widehat{\mu}(\xi)\widehat f(\xi)\widehat\psi(R^{-1}\xi)$ and both
$\widehat f$ and $\widehat\psi$ are bounded below on $\{|\xi|\lesssim R\}$, it follows that
$$
|\widehat g(\xi)|^{2}\gtrsim R^{-1}
$$
for at least $c R^{2}$ frequencies with $|\xi|\sim R$.
Moreover, for all $|\xi|\lesssim R$ we have the uniform bound
$$
|\widehat g(\xi)|^{2}\lesssim R^{-1} ,
$$
since $|\widehat{\mu}(\xi)|^{2}\lesssim R^{-1}$ for $|\xi|\sim R$ and $\widehat g$ vanishes for $|\xi|\gg R$.

Now let $P$ be a trigonometric polynomial of degree $k$, so that $\widehat P$ is supported on a set
$\Lambda$ of at most $k$ frequencies.
If
$$
\|g-P\|_{2}\le \eta\|g\|_{2}
$$
for some fixed $\eta\in (0,1)$, then almost all of the $L^{2}$ energy of $g$ must lie in $\Lambda$.
The total $L^{2}$ mass of $g$ satisfies
$$
\|g\|_{2}^{2} = \sum_{\xi} |\widehat g(\xi)|^{2}\approx R^{2}\cdot R^{-1}=R .
$$
For the contribution from the frequencies in $\Lambda$, we use the uniform bound on $|\widehat g|$ to obtain
$$
\sum_{\xi\in\Lambda} |\widehat g(\xi)|^{2}
\le \sum_{\xi\in\Lambda} C R^{-1}
\le k\cdot C R^{-1} .
$$
If $k=o(R^{2})$, then
$$
\frac{\sum_{\xi\in\Lambda} |\widehat g(\xi)|^{2}}{\|g\|_{2}^{2}}
\le \frac{k\cdot C R^{-1}}{c R}
= o(1) .
$$
Hence
$$
\frac{\sum_{\xi\notin\Lambda} |\widehat g(\xi)|^{2}}{\\|g\|_{2}^{2}}
= 1-o(1) ,
$$
and therefore
$$
\|g-P\|_{2}\ge (1-o(1))^{1/2}\,\|g\|_{2} .
$$
This contradicts $\|g-P\|_{2}\le \eta\|g\|_{2}$ for any fixed $\eta<1$ once $R$ is large.
Thus any such approximation requires $k\ge c(\eta) R^{2-o(1)}$. This completes the proof.

\vskip.125in 

\subsection{Proof of Theorem \ref{theorem:convexstationaryphase}} 

We write the Fourier transform in local coordinates.  Since $\partial K$ is convex, it may be covered by finitely many coordinate charts in which $\partial K$ is the graph of a convex function.  In each chart we may parametrize $\partial K$ by a Lipschitz map $\Gamma:U\subset{\mathbb R}^{d-1}\to{\mathbb R}^d$ and write
$$
\widehat{\mu}(\xi)=\int_U e^{-2\pi i \Gamma(u)\cdot \xi} J(u)\, du,
$$
where $J(u)$ is the surface Jacobian, bounded and measurable.

Fix $\xi$ with $|\xi|\sim R$ and $\xi/|\xi|\notin N(K)$.  
Since $\xi$ is separated from all outer normals, there exists a constant $c>0$ and an index $j\in\{1,\dots,d-1\}$ such that
$$
|\partial_{u_j}(\Gamma(u)\cdot \xi)| \ge c|\xi|
$$
for all $u\in U$.  
(This uses only convexity: the angle between $\xi$ and every outer normal is uniformly positive on a compact covering of $\partial K$.)

We now mollify $\Gamma$.  Let $\Gamma_\varepsilon=\Gamma*\rho_\varepsilon$, where $\rho_\varepsilon$ is a standard mollifier.  
Then $\Gamma_\varepsilon$ is smooth, $\Gamma_\varepsilon\to\Gamma$ uniformly, and  
$\partial_{u_j}\Gamma_\varepsilon\to\partial_{u_j}\Gamma$ in $L^\infty$ since $\Gamma$ is Lipschitz.  
Define
$$
I_\varepsilon(\xi)=\int_U e^{-2\pi i \Gamma_\varepsilon(u)\cdot\xi} J(u)\,du.
$$

Because $\partial_{u_j}(\Gamma_\varepsilon(u)\cdot\xi)$ converges uniformly to $\partial_{u_j}(\Gamma(u)\cdot\xi)$ and the latter is bounded below by $c|\xi|$, we have
$$
|\partial_{u_j}(\Gamma_\varepsilon(u)\cdot\xi)|\ge c|\xi|/2
$$
for all sufficiently small $\varepsilon$.

Since $\Gamma_\varepsilon$ is smooth, we may integrate by parts repeatedly:
$$
I_\varepsilon(\xi)
=\int_U \Big(L_\varepsilon^N e^{-2\pi i \Gamma_\varepsilon(u)\cdot\xi}\Big) J(u)\,du,
$$
where
$$
L_\varepsilon=\frac{1}{-2\pi i\,\partial_{u_j}(\Gamma_\varepsilon(u)\cdot\xi)}\,\partial_{u_j}.
$$
Each application of $L_\varepsilon$ contributes a factor bounded by $C|\xi|^{-1}$, so
$$
|I_\varepsilon(\xi)| \le C_N |\xi|^{-N}.
$$

We now pass to the limit $\varepsilon\to0$.  
Since $\Gamma_\varepsilon\to\Gamma$ uniformly and $J$ is bounded, dominated convergence gives
$$
I_\varepsilon(\xi)\to \widehat{\mu}(\xi).
$$
Thus
$$
|\widehat{\mu}(\xi)| \le C_N |\xi|^{-N}.
$$

See \cite{Str83} (Lemma 2.2) for a similar argument. 

\vskip.125in 

Integrating over the annulus $\{|\xi|\sim R\}$ gives
$$
\int_{\{|\xi|\sim R,\ \xi/|\xi|\notin N(K)\}} |\widehat{\mu}(\xi)|\, d\xi
\le C_N R^{d-1} R^{-N}=C_N R^{-(N-d+1)}.
$$
Since $N$ is arbitrary, this proves the first claim.  

For the second claim, since $N(K)\subset N(K)_{R^{-1}}$, we have $X_R^c \subset \{\xi:|\xi|\sim R,\ \xi/|\xi|\notin N(K)\}$. Hence the above decay estimate controls the entire complement of $X_R$.  

Moreover, $\|\widehat{\mu}_{R^{-1}}\|_{L^1}$ is $\gtrsim R^{(d-1)/2}$ (coming from the contribution of any smooth point of $\partial K$).  
Thus the ratio
$$
\frac{\|\widehat{\mu}_{R^{-1}}\|_{L^1(X_R^c)}}{\|\widehat{\mu}_{R^{-1}}\|_{L^1}}
$$
tends to zero as $R\to\infty$, and the desired $L^1$ concentration \eqref{eq:convexL1concentration_modified} follows.
\vskip.125in 

\subsection{Proof of Corollary \ref{cor:minkowskikicksass}} 

By Theorem \ref{theorem:fractalsandwich}, 
$$ \FR(f\mu) \leq 2 \cdot \sqrt{\frac{|X_R|}{R^d}}, $$ where we set $\eta=\frac{1}{2}$ in Theorem \ref{theorem:convexstationaryphase}. 

By assumption and the definition of upper Minkowski dimension, for any $\epsilon>0$ there exists $C_{\epsilon}>0$ such that 
$$ |X_R| \leq C_{\epsilon} R^{a+1+\epsilon},$$ where $a$ is the upper Minkowski dimension of $N(K)$. It follows that 
$$ \FR(f\mu) \leq 2 \cdot C_{\epsilon}^{\frac{1}{2}} \cdot \sqrt{R^{a-(d-1)+\epsilon}}. $$

\vskip.125in 

By Theorem \ref{theorem:L2approximation}, there exists a trigonometric polynomial of degree 
$$ \eta^{-2} \cdot C_{\epsilon} \cdot R^{a-(d-1)+\epsilon} \cdot R^d=\eta^{-2} \cdot C_{\epsilon} \cdot R^{a+1+\epsilon},$$ such that 
$$ {\left|\left|{(f\mu)}_{R^{-1}}-P\right|\right|}_2 \leq \eta \cdot {\left|\left|{(f\mu)}_{R^{-1}}\right|\right|}_2.$$

\vskip.125in 

This completes the proof.


\newpage 

\begin{thebibliography}{8}

\bibitem{A2025} K. Aldaleh, W. Burstein, G. Garza, A. Iosevich, J. Iosevich, A. Khalil, J. King, T. Le, I. Li, A. Mayeli, K. Nguyen, and N. Shaffer, {\it The Fourier Ratio and complexity of signals}, (in preparation), (2025). 

\bibitem{Bourgain89} J. Bourgain, {\it Bounded orthogonal systems and the $\Lambda(p)$-set problem}, Acta Math. \textbf{162} (1989), no. 3-4, 227–245.

\bibitem{BD18} J. Bourgain and S. Dyatlov, {\it Spectral gaps without the pressure condition}, Ann. of Math. (2) \textbf{187} (2018), no. 3, 825-867.

\bibitem{BCT97} L. Brandolini, L. Colzani, and G. Travaglini, {\it Average decay of Fourier transforms and integer points in polyhedra}, Arkiv fur Matematik 35 (1997), no. 2, 253–275.

\bibitem{BIT03} L. Brandolini, A. Iosevich, and G. Travaglini, {\it Planar convex bodies, Fourier transform, lattice points, and irregularities of distribution}, Transactions of the American Mathematical Society 355 (2003), no. 9, 3513–3535.

\bibitem{BRT98} L.~Brandolini, M.~Rigoli and G.~Travaglini, {\it Average decay of Fourier transforms and geometry of convex sets}, Rev. Mat. Iberoamericana 14 (1998), no. 3, 519–560.  

\bibitem{BIMN2025} W. Burstein, A. Iosevich, A. Mayeli, and H. Nathan, {\it Fourier minimization and time series imputation}, (arXiv:2506.19226), (2025). 

\bibitem{DS89} D. Donoho and P. Stark, {\it Uncertainty principle and signal processing}, SIAM Journal of Applied Math., (1989), Society for Industrial and Applied Mathematics, volume 49, No. 3, pp. 906-931. 

\bibitem{Hagerstrom2025} W. Hagerstrom, {\it A number of perspectives on signal recovery}, University of Rochester Honors Thesis (2025). 

\bibitem{LW18} I. Laba and H. Wang, {\it Decoupling and near-optimal restriction estimates for Cantor sets}, (English summary)
Int. Math. Res. Not. IMRN (2018), no. 9, 2944–2966.

\bibitem{Stein1993} E.M. Stein, {\it Harmonic Analysis: Real-Variable Methods, Orthogonality, and Oscillatory Integrals}, Princeton Mathematical Series, vol. 43, Princeton Univ. Press, Princeton, NJ, 1993.

\bibitem{Str83}
R. S. Strichartz, 
Restrictions of Fourier transforms to quadratic surfaces and decay of solutions of wave equations, 
Duke Math. J. 50 (1983), 549–567.

\bibitem{Strichartz83} R.S. Strichartz, {\it Convolutions with kernels having singularities on a sphere},
Trans. Amer. Math. Soc. \textbf{148} (1970), 461-471.

\bibitem{Talagrand98} M. Talagrand, {\it Selecting a proportion of characters}, Israel J. Math. \textbf{108} (1998),
173-191.

\bibitem{Tomas1975} P. A. Tomas, {\it A restriction theorem for the Fourier transform}, Bull. Amer. Math. Soc. \textbf{81} (1975), no.~2, 477-478.

\end{thebibliography}
\end{document}